\documentclass[11pt,a4paper,oneside,article]{article}
\usepackage[T2A]{fontenc} 
\usepackage[utf8]{inputenc} \usepackage{amsfonts,amsmath, amsthm} 
\usepackage{amssymb} 
\usepackage{algorithm}
\usepackage{algorithmic} 
\usepackage{xcolor} 
\usepackage{url}
\usepackage{etoolbox} 
\usepackage{pdfpages}
\usepackage{hyperref}

\usepackage{natbib}

\newcommand{\R}{{\mathbb R}}
\newcommand{\E}{{\mathbb E}}
\newcommand{\Prob}{{\mathbb P}}
\def\e{\varepsilon}
\newcommand{\ee}{{\rm e}}
\newcommand{\la}{\langle}
\newcommand{\ra}{\rangle}

\def\vp{\varphi}
\def\tpi{\tilde{\pi}}
\def\tPi{\tilde{\Pi}}
\def\S{{\mathcal S}}
\def\U{\mathcal{U}}
\def\Ec{\mathcal E}

\def\EE{\mathbb{E}}

\def\tf{\tilde{f}}
\def\tg{\tilde{g}}

\def\tf{\widetilde{f}}

\def\tnmf{\widetilde{\nabla}^m f}
\newcommand{\argmin}{\mathop{\arg\!\min}}
\def\O{\widetilde{O}}

\newcommand{\V}{\mathcal{V}}
\newcommand{\B}{{\mathcal B}}
\def\C{{\mathcal C}}
\newcommand{\Bb}{{\mathfrak B}}
\newcommand{\diam}{{\rm diam}}

\def\one{{\mathbf 1}}
\newcommand{\transp}[1]{1^T}
\def\Reg{\mathcal R}

\def\W{\mathcal W}
\def\M{\mathcal M_+^1}
\def\X{\mathcal X}
\def\Y{\mathcal Y}

\def\Blm{\boldsymbol{\lambda}}

\def\Bzeta{\boldsymbol{\zeta}}
\def\Beta{\boldsymbol{\eta}}
\def\BBlm{\boldsymbol{\bar{\lambda}}}

\def\BBzeta{\boldsymbol{\bar{\zeta}}}
\def\BBeta{\boldsymbol{\bar{\eta}}}
\def\p{\mathtt{p}}
\def\blm{\bar{\lambda}}

\def\eig{\lambda_{\max}(W)}

\newcommand{\bm}[1]{\mathbf{#1}}
\newcommand{\eps}{\varepsilon}

\newtheorem{Lm}{Lemma}[section]

\newtheorem{Th}{Theorem}[section]

\providetoggle{thesis} \settoggle{thesis}{false}

\begin{document}

\begin{titlepage}
\iftoggle{thesis}
{
   \begin{center}
Moscow Institute of Physics and Technology
       \vspace{0.5cm}
       \begin{flushright}
       \it
as a manuscript
       \end{flushright}
       
       \vspace{1.5cm}
       
       Pavel Evgenievich Dvurechenskii
       
\vspace{1.5cm}       
\textbf{\Large Numerical methods in large-scale optimization: inexact oracle and primal-dual analysis}
       
       \vspace{4.5cm}
       
       Dissertation
       
       for the purpose of obtaining academic degree
       
       Doctor of Science in Computer Science

       \vfill
            
       Moscow - 2020
            
   \end{center}
   }
   {
   \begin{center}
Moscow Institute of Physics and Technology
       \vspace{0.5cm}
       \begin{flushright}
       \it
as a manuscript
       \end{flushright}
       
       \vspace{1.5cm}
       
       Pavel Evgenievich Dvurechenskii
       
\vspace{1.5cm}       
\textbf{\Large Numerical methods in large-scale optimization: inexact oracle and primal-dual analysis}
       
       \vspace{4.5cm}
       
       Dissertation Summary
       
       for the purpose of obtaining academic degree
       
       Doctor of Sciences in Computer Science

       \vfill
            
       Moscow - 2020
            
   \end{center}
   }

\end{titlepage}

\iftoggle{thesis}
{
{\small\tableofcontents}
\newpage
}
{
The dissertation was prepared at Moscow Institute of Physics and Technology.
\\

{\bf Scientific Consultant:} \\

Alexander Vladimirovich Gasnikov, Doctor of Sciences in Mathematical Modelling, Numerical Methods and Software Complexes, 
Associate Professor at Mathematical Foundations of Control chair, Moscow Institute of Physics and Technology.

\newpage
}

\newpage
\section{Introduction}

Numerical optimization remains an active area of research since 1980's,  motivated by a vast range of applications, e.g. operations research, optimal control. Starting with the works \cite{karmarkar1984new,nesterov1994interior} one of the main areas of research in numerical optimization became interior-point methods. These methods combine Newton steps with penalty approach and allow to solve a very general class of convex problems in polynomial-time, which is justified both theoretically and practically. The new century introduced new challenges for numerical methods in optimization. Thanks to increasing amount of available data and more powerful computational resources, machine learning became an area of intensive research. A cornerstone optimization problem in machine learning is the empirical risk minimization with the key aspect being large dimension of the decision variable and large number of components used in the objective function. In this setting the Newton iteration becomes expensive in general since it requires matrix inversion. This motivated a sacrifice of logarithmic dependence on the accuracy to a cheap iteration and the use of first-order methods to solve such problems. Another reason was that the data is usually noisy and there is no need to solve the optimization problem to a high accuracy in this setting. Another main application for first-order methods is signal processing and image analysis, where the goal is to reconstruct a high-dimensional signal from high-dimensional data, e.g. noisy images.

Yet, known already for a long time \cite{cauchy1847methode,polyak1963gradient,robbins1951stochastic}, first-order methods entered their renaissance in 2000's. Some important facts on these methods were already known for 15 years. In particular, the concept of black-box oracle \cite{nemirovsky1983problem} allowed to obtain lower worst-case complexity bounds for different classes of problems and methods. In particular, a gap was discovered between the lower bound $O(1/k^2)$ and an upper bound  $O(1/k)$ for gradient method for minimizing convex smooth functions. Here $k$ is the iteration counter. This gap led to an important phenomenon of acceleration for first-order methods and accelerated gradient method \cite{nesterov1983method}. In the new century many extensions of this algorithm were proposed motivated by image processing problems and machine learning, including composite versions \cite{beck2009fast,nesterov2013gradient}, accelerated stochastic gradient method \cite{lan2012optimal}, accelerated variance reduction methods \cite{johnson2013accelerating,lin2014accelerated,lin2015universal,lan2017optimal,shalev-shwartz2014accelerated}. In addition to accelerated stochastic gradient methods for finite-sum problems, which use a random choice of the gradient of the component, acceleration was introduced for other randomized methods such as random coordinate descent \cite{nesterov2012efficiency} and random gradient-free optimization \cite{nesterov2017random}. The latter is motivated by problems, in which only zero-order oracle is available, e.g. when the objective is given as a solution of some auxiliary problem. For this setting, it is important to analyze zero-order methods with inexact function values since this auxiliary problem may be possible to solve only inexactly. In the setting of first-order methods inexactness may also be encountered in practice. Accelerated gradient method with inexact gradients was analyzed in \cite{aspremont2008smooth}, and an important framework of inexact first-order oracle was introduced in \cite{devolder2014first} and several extensions are discussed in \cite{gasnikov2016stochastic}. Another important extension of accelerated first-order methods are accelerated methods for problems with linear constraints, which was proposed in \cite{beck2014fast}, yet with a non-optimal rate $O(1/k)$ for the constraints feasibility.

{\bf  Object and goals of the dissertation.} The goal of the dissertation is twofold. The first goal is to further extend the existing first and zero-order methods for problems with inexactness in function and gradient values, the inexactness being deterministic or stochastic. The second goal is to construct new primal-dual first-order methods, which allow to solve simultaneously the primal and dual problem with optimal convergence rates. A particular focus is made on problems with linear constraints and the application of the proposed methods to optimal transport distance and barycenter problems.  

{\bf  The obtained results:}
\begin{enumerate}
    \item We propose a stochastic intermediate gradient method for convex problems with stochastic inexact oracle.
    \item We develop a gradient method with inexact oracle for deterministic non-convex optimization.
    \item We develop gradient-free method with inexact oracle for deterministic convex optimization.
    \item We develop a method to calculate the derivative of the pagerank vector and in combination with the above two methods propose gradient-based and gradient-free optimization methods for learning supervised pagerank model.
    \item We develop a concept of inexact oracle for the methods which use directional derivatives and propose accelerated directional derivative method for smooth stochastic convex optimization. We also develop an accelerated and non-accelerated directional derivative method for strongly convex smooth stochastic optimization.
    \item We develop primal-dual methods for solving infinite-dimensional games in convex-concave and strongly convex-concave setting.
    \item We develop non-adaptive and adaptive accelerated primal-dual gradient method for strongly convex minimization problems with linear equality and inequality constraints.
    \item We apply this algorithm to the optimal transport problem and obtain new complexity estimates for this problem, which in some regime are better than the ones for the Sinkhorn's algorithm.
    \item We propose a stochastic primal-dual accelerated gradient method for problems with linear constraints and apply it to the problem of approximation of Wasserstein barycenter.
    \item We propose a primal-dual extension of accelerated methods which use line-search to define the stepsize and to be adaptive to the Lipschitz constant of the gradient.
\end{enumerate}

{\bf Author’s contribution} includes the development of the listed above optimization methods, proving convergence rates and complexity result theorems for these methods and their applications to optimal transport problems and learning problem for a supervised pagerank model.

{\bf Novelties.}
The proposed versions of accelerated  first and zero-order methods for convex optimization under different types of inexactness are novel. The proposed primal-dual methods for the listed setups are also novel, and allow to obtain new methods for optimal transport problems. In particular, we obtain new complexity results for non-regularized optimal transport problem and a new distributed algorithm for approximating Wasserstein barycenter of a set of measures using samples from these measures.

As a result of the work on this dissertation, 10 papers were published:

{\bf First-tier publications:}

\begin{enumerate}
    \item Dvurechensky, P., and Gasnikov, A. Stochastic intermediate gradient method for convex problems  with stochastic inexact oracle. Journal of Optimization Theory and Applications 171, 1 (2016), 121--145, Scopus Q1 (main co-author; the author of this thesis proposed main algorithms, formulated and proved convergence rate theorems for the proposed methods).
    \item Gasnikov, A. V., and Dvurechensky, P. E. Stochastic intermediate gradient method for convex  optimization problems. Doklady Mathematics 93, 2 (2016), 148--151, Scopus Q2 (main co-author; the author of this thesis proposed main algorithms, formulated and proved convergence rate theorems for the proposed methods).
    \item Bogolubsky, L., Dvurechensky, P., Gasnikov, A., Gusev, G., Nesterov, Y., Raigorodskii, A. M., Tikhonov, A., and Zhukovskii, M. Learning supervised pagerank with gradient-based and  gradient-free optimization methods. In Advances in Neural Information Processing Systems 29, D. D. Lee, M.  Sugiyama, U. V. Luxburg, I. Guyon, and R. Garnett, Eds. Curran Associates, Inc., 2016, pp. 4914--4922, CORE A* (the author of this thesis proposed general gradient-free (Algorithm 1,2) and gradient (Algorithm 3,4) methods with inexact oracle, proposed a method for approximating the derivative of the pagerank vector, formulated and proved convergence rate theorems for the proposed methods: Lemma 1,2, Theorem 1-4).
    \item Dvurechensky, P., Gorbunov, E., and Gasnikov, A. An accelerated directional derivative method  for smooth stochastic convex optimization. European Journal of Operational Research (2020), \url{https://doi.org/10.1016/j.ejor.2020.08.027}, Scopus Q1 (main co-author; the author of this thesis proposed a concept of inexact oracle for directional derivatives in stochastic convex optimization, proved (in inseparable cooperation with E. Gorbunov) convergence rate Theorem 1 for the accelerated directional derivative method, proved convergence rate Theorems 3,4 for strongly convex problems).
    \item Dvurechensky, P., Nesterov, Y., and Spokoiny, V. Primal-dual methods for solving in infinite-dimensional games. Journal of Optimization Theory and Applications 166, 1 (2015), 23--51, Scopus Q1 (main co-author; the author of this thesis developed main algorithms and proved convergence rate theorems).
    \item Dvurechensky, P., Gasnikov, A., and Kroshnin, A. Computational optimal transport: Complexity by accelerated gradient descent is better than by Sinkhorn's algorithm. In Proceedings of the  5th  International Conference on Machine Learning (2018), J. Dy and A. Krause, Eds., vol. 80 of Proceedings of  Machine Learning Research, pp. 1367--1376, CORE A* (main co-author; the author of this thesis proposed general primal-dual adaptive accelerated gradient method (Algorithm 3) for problems with linear constraints, proved convergence rate Theorem 3, proposed an algorithm for approximating optimal transport (OT) distance (Algorithm 4), obtained complexity bound for approximating OT distance (Theorem 4), performed numerical experiments for comparison of this method with the Sinkhorn's method).
    \item Dvurechensky, P., Dvinskikh, D., Gasnikov, A., Uribe, C. A., and Nedi\'c, A. Decentralize and randomize: Faster algorithm for Wasserstein barycenters. In Advances in Neural Information  Processing Systems 31 (2018), S. Bengio, H. Wallach, H. Larochelle, K. Grauman, N. Cesa-Bianchi, and R. Garnett,  Eds., NeurIPS 2018, Curran Associates, Inc., pp. 10783--10793, CORE A* (main co-author; the author of this thesis proposed the general idea of the paper, general primal-dual accelerated stochastic gradient method (Algorithm 2) for problems with linear constraints, proved convergence rate Theorem 2, proposed an algorithm for approximating Wasserstein barycenter (Algorithm 4), proved (in inseparable cooperation with D. Dvinskikh) its complexity Theorem 3).
    \item Guminov, S. V., Nesterov, Y. E., Dvurechensky, P. E., and Gasnikov, A. V. Accelerated primal-dual gradient descent with linesearch for convex, nonconvex, and nonsmooth optimization  problems. Doklady Mathematics 99, 2 (2019), 125-128, Scopus Q2 (the author of this thesis proposed a primal-dual variant of the accelerated gradient method with linesearch for problems with linear constraints, proved convergence rate Theorem 3).
    \item Nesterov, Y., Gasnikov, A., Guminov, S., and Dvurechensky, P. Primal-dual accelerated gradient  methods with small-dimensional relaxation oracle. Optimization Methods and Software (2020), \url{https://doi.org/10.1080/10556788.2020.1731747}, Scopus Q1 (the author of this thesis proposed a primal-dual variant of the universal accelerated gradient method with small-dimensional relaxation (Algorithm 7) for problems with linear constraints, proved its convergence rate Theorem 4.1). 
\end{enumerate}

{\bf Second-tier publications:}

\begin{enumerate}
   \item Chernov, A., Dvurechensky, P., and Gasnikov, A. Fast primal-dual gradient method for strongly  convex minimization problems with linear constraints. In Discrete Optimization and Operations Research:  9th International Conference, DOOR 2016, Vladivostok, Russia, September 19-23, 2016, Proceedings (2016),  Y. Kochetov, M. Khachay, V. Beresnev, E. Nurminski, and P. Pardalos, Eds., Springer International Publishing,  pp. 391--403, Web of Science and Scopus (main co-author; the author of this thesis developed main algorithm and proved convergence rate theorem).
\end{enumerate}

{\bf Reports at conferences and seminars:}

\begin{enumerate}
    \item International Workshop "Advances in Optimization and Statistics", Berlin, 15.05.2014--16.05.2014, "Stochastic Intermediate Gradient Method for Convex Problems with Inexact Stochastic Oracle".
    \item Seminar "Modern Methods in Applied Stochastics and Nonparametric Statistics", Berlin, 03.06.2014, "Gradient methods for convex problems with stochastic inexact oracle".
    \item V International Conference on Optimization Methods and Applications (OPTIMA-2014), Petrovac, Montenegro, 28.09.2014--04.10.2014, "Gradient-free optimization methods with ball randomization".
    \item  VI traditional school for young scientists "Control, information, optimization", Moscow, 22.06.2014-29.06.2014, "Gradient methods for convex problems with stochastic inexact oracle".
    \item 38-th conference-school of IITP RAS "Information technologies and systems", Nizhnii Novgorod, 01.09.2014--05.09.2014, "Stochastic Intermediate Gradient Method for Convex Problems with Inexact Stochastic Oracle".
    \item Workshop “Frontiers of High Dimensional Statistics, Optimization, and Econometrics”, Moscow, 26.02.2015--27.02.2015 г., "Random gradient-free methods for random walk based web page ranking functions learning".
    \item VII traditional school for young scientists "Control, information, optimization", Moscow, 14.06.2014-20.06.2014, "Semi-Supervised PageRank Model Learning with Gradient-Free Optimization Methods".
    \item 29-th conference-school of IITP RAS "Information technologies and systems", Sochi, 07.09.2014--11.09.2015, "Stochastic Intermediate Gradient Method: convex and strongly-convex case".
    \item 30th annual conference of Belgian Operational Research Society (ORBEL 30), Louvain-la-Neuve, Belgium, 28.01.2016--29.01.2016, "Random gradient-free methods for ranking algorithm learning".
    \item Workshop on Modern Statistics and Optimization, Moscow, 23.02.2016--24.02.2016, "Gradient and gradient-free methods for pagerank algorithm learning".
    \item VII International Conference Optimization and Applications (OPTIMA 2016), Petrovac, Montenegro, 25.09.2016--02.10.2016, "Accelerated Primal-Dual Gradient Method for Linearly Constrained Minimization Problems".
    \item VIII Moscow International Conference on Operations Research (ORM 2016), Moscow, 17.10.2016--22.10.2016, "Accelerated Primal-Dual Gradient Method for Composite Optimization with Unknown Smoothness Parameter"
    \item \textbf{Conference on Neural Information Processing Systems (NIPS 2016)}, Barcelona, 05.12.2016--10.12.2016, "Learning Supervised PageRank with Gradient-Based and Gradient-Free Optimization Methods".
    \item Workshop Shape, Images and Optimization, M\"unster, Germany, 28.02.2017--03.03.2017, "Gradient Method With Inexact Oracle for Composite Non-Convex Optimization".
    \item Optimization and Statistical Learning, Les Houches, France, 10.04.2017-14.04.2017, "Gradient Method With Inexact Oracle for Composite Non-Convex Optimization".
    \item Foundations of Computational Mathematics, Barcelona, Spain, 10.07.2017-19.07.2017, "Gradient Method With Inexact Oracle for Composite Non-Convex Optimization". 
    \item Co-Evolution of Nature and Society Modelling, Problems \& Experience. Devoted to Academician Nikita Moiseev centenary (Moiseev-100), Moscow, 07.11.2017--10.11.2017, "Adaptive Similar Triangles Method: a Stable Alternative to Sinkhorn’s Algorithm for Regularized Optimal Transport".
    \item 18th French-German-Italian Conference on Optimization, Germany, 25.09.2017--28.09.2017, Paderborn, Germany, "Gradient method with inexact oracle for composite non-convex optimization"
    \item 3. International Matheon Conference on Compressed Sensing and its Applications, Berlin, 04.12.2017--08.12.2017, "Adaptive Similar Triangles Method:  a Stable Alternative to Sinkhorn’s Algorithm for Regularized Optimal Transport".
    \item Games, Dynamics and Optimization (GDO2018), Vienna, Austria, 13.03.2018--15.03.2018, "Primal-Dual Methods for Solving Infinite -Dimensional Games".
    \item \textbf{International Conference on Machine Learning (ICML 2018)}, Stockholm, Sweden, 10.07.2018--15.07.2018, "Computational optimal transport: Complexity by accelerated gradient descent is better than by Sinkhorn's algorithm".
    \item 23rd International Symposium on Mathematical Programming, 01.07.2018--06.07.2018, "Computational Optimal Transport: Accelerated Gradient Descent vs Sinkhorn".
    \item Grenoble Optimization Days 2018: Optimization algorithms and applications in statistical learning, Grenoble, France, 28.06.2018--29.06.2018, "Faster algorithms for (regularized) optimal transport". 
    \item Statistical Optimal Transport Conference, Moscow, 24.07.2018--25.07.2018, "Computational Optimal Transport: Accelerated Gradient Descent vs Sinkhorn’s Algorithm".
    \item \textbf{Conference on Neural Information Processing Systems (NIPS 2018)}, Montreal, Canada, 02.12.2018--08.12.2018, "Decentralize and randomize: Faster algorithm for Wasserstein barycenters".
    \item Optimization and Statistical Learning, Les Houches, France, 24.03.2019--29.03.2019, "Distributed optimization for Wasserstein barycenter".
    \item \textbf{International Conference on Machine Learning (ICML 2019)}, Long Beach, USA, 09.06.2019--15.06.2019, "On the Complexity of Approximating Wasserstein Barycenters".
    \item International Conference on Continuous Optimization (ICCOPT 2019), Berlin, Germany, 03.08.2019-08.08.2019, "A Unifying Framework for Accelerated Randomized Optimization Methods".
    \item Workshop on optimization and applications, Moscow, 27.09.2019, "Accelerated Alternating Minimization for Optimal Transport".
    \item Recent advances in mass transportation, Moscow, 23.09.2019--27.09.2019, "On the complexity of optimal transport problems".
    \item Workshop by the GAMM Activity Group on Computational and Mathematical Methods in Data Science, Berlin, Germany, 24.10.2019--25.10.2019, "On the complexity of optimal transport problems".
    \item HSE-Yandex autumn school on generative models, Moscow, 26.11.2019--29.11.2019, "Optimization methods for optimal transport".
    \item Workshop on Mathematics of Deep Learning 2019, Berlin, Germany, 03.12.2019--05.12.2019, "On the complexity of optimal transport problems".
    \item Workshop on PDE Constrained Optimization under Uncertainty and Mean Field Games, Berlin, Germany, 28.01.2020--30.01.2020, "Distributed optimization for Wasserstein barycenters".
\end{enumerate}

\section{Optimization with inexact oracle}
In this section we briefly describe the methods and their convergence properties for optimization problems under inexact information. We consider first-order methods and directional derivative methods.

\subsection{Stochastic intermediate gradient method for convex problems  with stochastic inexact oracle}
The results of this subsection are published in \cite{dvurechensky2016stochastic,gasnikov2016stochasticInter}.

Let $E$ be a finite-dimensional real vector space and $E^*$ be its dual. We denote the value of a linear function $g \in E^*$ at $x\in E$ by $\la g, x \ra$. Let $\|\cdot\|$ be some norm on $E$. We denote by $\|\cdot\|_*$ the dual norm for $\|\cdot\|_E$, i.e. \\
$
\|g\|_*~=~\sup_{y \in E} \{ \la g,y \ra:~\|y\|_E\leq 1 \}
$. By $\partial f(x)$ we denote the subdifferential of the function $f(x)$ at a point $x$.
In this subsection, we consider the {\it composite optimization} problem of the form
\begin{equation}
\min_{x \in Q} \{ \varphi(x) := f(x) + h(x)\},
\label{eq1:PrStateInit} 
\end{equation}
where $Q \subset E$ is a closed and convex set, $h(x)$ is a simple convex function, $f(x)$ is a convex function with { \it stochastic inexact oracle } \cite{devolder2011stochastic}. This means that, for every $x \in Q$, there exist $f_{\delta,L}(x) \in \R$ and $g_{\delta,L}(x) \in E^*$, such that
\begin{equation}
0 \leq f(y) - f_{\delta,L}(x) - \la g_{\delta,L}(x) , y-x \ra \leq \frac{L}{2} \|x-y\|^2 + \delta, \quad \forall y \in Q,
\label{eq:dLOracle}
\end{equation}
and also that, instead of $(f_{\delta,L}(x), g_{\delta,L}(x))$ (we will call this pair a $(\delta,L)$-oracle), we use their stochastic approximations $(F_{\delta,L}(x, \xi), G_{\delta,L}(x, \xi))$. 
The latter means that, for any point $x \in Q$, we associate
with $x$ a random variable $\xi$ whose probability distribution is supported on a set $\Xi \subset \R$ and such that $\E_{\xi} F_{\delta,L}(x, \xi) = f_{\delta,L}(x)$, $\E_{\xi} G_{\delta,L}(x, \xi) = g_{\delta,L}(x)$ and
$
 \E_{\xi} (\|G_{\delta,L}(x, \xi) - g_{\delta,L}(x)\|_*)^2 \leq \sigma^2. $

To deal with such problems we will need a {\it prox-function} $d(x)$, which is differentiable and strongly convex with parameter 1 on $Q$ with respect to $\|\cdot\|$. Let $x_0$ be the minimizer of $d(x)$ on $Q$. By translating and scaling $d(x)$, if necessary, we can always ensure that $
d(x_0) = 0$, $d(x) \geq \frac12 \|x-x_0\|^2$, $\forall x \in Q$.
We define also the corresponding {\it Bregman distance}:
$V (x, z) = d(x) - d(z) - \la \nabla d(z), x - z \ra$.
Let $\{\alpha_i\}_{i \geq 0}$, $\{\beta_i\}_{i \geq 0}$, $\{B_i\}_{i \geq 0} \subset \R$ be three sequences of coefficients satisfying 
\begin{align}
& \alpha_0 \in ]0,1], \quad \beta_{i+1} \geq \beta_i > L , \quad \forall i \geq 0, \label{eq:alpbet1} \\
& 0\leq \alpha_i \leq B_i , \quad \forall i \geq 0, \label{eq:alpB} \\
& \alpha_k^2 \beta _k  \leq B_k \beta_{k-1} \leq \left( \sum_{i=0}^k{\alpha_i} \right) \beta_{k-1}, \quad \forall k \geq 1.
 \label{eq:alpbetB} \\
 & A_k := \sum_{i=0}^k{\alpha_i}, \;\;\;\tau_{i}:=\frac{\alpha_{i+1}}{B_{i+1}}
\end{align}
The Stochastic Intermediate Gradient Method (SIGM) is described below as Algorithm \ref{alg:SIGM}.  Let $a \geq 1$ and $b \geq 0$ be some parameters. Let us assume that we know a number $R$ such that $ \sqrt{2d(x^*)} \leq R$. We set for $p\in [1,2]$
\begin{align}
& \alpha_i = \frac{1}{a} \left(\frac{i+p}{p} \right)^{p-1}, \quad  \forall i \geq 0, \label{eq:chooseAlp}\\
& \beta_i = L + \frac{b\sigma}{R} (i+p+1)^{\frac{2p-1}{2}} , \quad  \forall i \geq 0, \label{eq:chooseBet}\\
& B_i = a \alpha_i^2 = \frac{1}{a} \left(\frac{i+p}{p} \right)^{2p-2}, \quad  \forall i \geq 0. \label{eq:chooseB}
\end{align}

\begin{Th}
If the sequences $\{\alpha_i\}_{i \geq 0}$, $\{\beta_i\}_{i \geq 0}$, $\{B_i\}_{i \geq 0}$ are chosen according to \eqref{eq:chooseAlp}, \eqref{eq:chooseBet}, \eqref{eq:chooseB} with $a = 2^{\frac{2p-1}{2}}$ and $b=2^{\frac{5-2p}{4}}p^{\frac{1-2p}{2}}$, then the
sequence $y_k$ generated by the SIGM satisfies
\begin{align}
& \E_{\xi_0, \dots, \xi_k} \varphi(y_k) - \varphi^* \leq  \frac{LR^2 p^p 2^{\frac{2p-3}{2}} }{(k+p)^{p}} + \frac{\sigma R 2^{\frac{3+2p}{4}  }\sqrt{p}(k+p+2)^{p-\frac12}}{(k+p)^{p}} +  \notag \\
& + 2^{2p-1} \left( \left(\frac{k+p}{p} \right)^{p-1} + 1 \right) \delta \leq \frac{C_1 LR^2  }{k^{p}} + \frac{C_2 \sigma R }{\sqrt{k}} +  C_3 k^{p-1}\delta =  \notag \\
& = \Theta\left(\frac{LR^2}{k^p} + \frac{\sigma R}{\sqrt{k}} + k^{p-1}\delta\right) ,\notag
\end{align}
where $C_1=4 \sqrt{2}$, $C_2 = 16 \sqrt{2}$, $C_3 = 48$.
\label{Th04}
\end{Th}

\begin{algorithm}[h!]
\caption{Stochastic Intermediate Gradient Method (SIGM)}
\label{alg:SIGM}
\begin{algorithmic}[1]
\REQUIRE {The sequences $\{\alpha_i\}_{i \geq 0}$, $\{\beta_i\}_{i \geq 0}$, $\{B_i\}_{i \geq 0}$, functions $d(x)$, $V(x,z)$.}
\ENSURE {The point $y_{k}$.}
\STATE  Compute $x_0 := \arg \min_{x\in Q} \{d(x)\}$. Let $\xi_0$ be a realization of the random variable $\xi$. Calculate $G_{\delta,L}(x_0,\xi_0)$. Set $k$ = 0.
\STATE 
$y_0 := \arg \min_{x\in Q} \{ \beta_0 d(x) + \alpha_0 \la G_{\delta,L}(x_0,\xi_0), x-x_0 \ra + \alpha_0h(x)\}$.
\REPEAT
\STATE $z_k := \arg \min_{x\in Q} \{ \beta_k d(x) + \sum_{i=0}^k{\alpha_i \la G_{\delta,L}(x_i,\xi_i), x-x_i \ra } + A_kh(x) \}$.
\STATE $x_{k+1} := \tau_k z_k + (1- \tau_k) y_k$.
\STATE Let $\xi_{k+1}$ be a realization of the random variable $\xi$. Calculate $G_{\delta,L}(x_{k+1},\xi_{k+1})$.
\STATE $\hat {x}_{k+1} := \arg \min_{x\in Q} \{ \beta_k V(x,z_k) + \alpha_{k+1} \la G_{\delta,L}(x_{k+1},\xi_{k+1}), x-z_k \ra + \alpha_{k+1} h(x).\}$.
\STATE $w_{k+1} := \tau_k \hat {x}_{k+1} + (1- \tau_k) y_k$.
\STATE $y_{k+1} := \frac{A_{k+1}-B_{k+1}}{A_{k+1}} y_{k} + \frac{B_{k+1}}{A_{k+1}} w_{k+1}$.
\UNTIL{}
\end{algorithmic}
\end{algorithm}

It is possible to obtain an upper bound on the probability of large deviations for the $\varphi(y_k) - \varphi^*$. To do that, we make the following additional assumptions.
\begin{enumerate}
	\item $\xi_0, \dots, \xi_k$ are i.i.d random variables. \label{as1}
	\item $G_{\delta,L}(x,\xi)$ satisfies the light-tail condition
	$$
	\E_\xi \left[ \exp \left( \frac{\|G_{\delta,L}(x,\xi)-g_{\delta,L}(x) \|^2_*}{\sigma^2} \right) \right] \leq \exp(1).
	$$ \label{as2}
	\item Set $Q$ is bounded, and we know a number $D>0$, such that $ \max_{x,y\in Q}\|x~-~y\| \leq D$. \label{as3}
\end{enumerate}

\begin{Th}
If the sequences $\{\alpha_i\}_{i \geq 0}$, $\{\beta_i\}_{i \geq 0}$, $\{B_i\}_{i \geq 0}$ are chosen according to \eqref{eq:chooseAlp}, \eqref{eq:chooseBet}, \eqref{eq:chooseB} with $a = 2^{\frac{2p-1}{2}}$ and $b=2^{\frac{5-2p}{4}}p^{\frac{1-2p}{2}}$, then the
sequence $y_k$ generated by the SIGM satisfies
\begin{align}
& \Prob \Biggl( \varphi(y_k) - \varphi^* >  \frac{C_1 LR^2 }{k^{p}} + \frac{C_2( 1+ \Omega )\sigma R }{\sqrt{k}}  + C_3 k^{p-1}  \delta + \frac{C_4 D \sigma \sqrt{\Omega} }{\sqrt{k}} \Biggr)   \notag \\
& \leq \Prob \Biggl( \varphi(y_k) - \varphi^* > \frac{LR^2 p^p 2^{\frac{2p-3}{2}} }{(k+p)^{p}} + \frac{( 1+ \Omega )\sigma R 2^{\frac{3+2p}{4}  }\sqrt{p}(k+p+2)^{p-\frac12}}{(k+p)^{p}}   \notag \\
& + 2^{2p-1} \left( \left(\frac{k+p}{p} \right)^{p-1} + 1 \right) \delta + \frac{2D \sigma \sqrt{6 \Omega p }}{\sqrt{k+p}} \Biggr) \notag  \leq 3 \exp(- \Omega)  , \notag
\end{align}
where $C_1=4 \sqrt{2}$, $C_2 = 16 \sqrt{2}$, $C_3 = 48$, $C_4 = 4 \sqrt{3}$.
\label{Th05}
\end{Th}

Next, we consider two modifications of the SIGM for strongly convex problems. For the first modification, we obtain the rate of convergence in terms of the non-optimality gap expectation and for the second we bound the probability of large deviations from this rate. 
We additionally assume that $E$ is a Euclidean space with scalar product $\la \cdot , \cdot \ra$ and norm $\|x\| := \sqrt{\la x , H x \ra}$, where $H$ is a symmetric positive definite matrix. Without loss of generality, we assume that the function $d(x)$ satisfies conditions $0=\arg \min_{x\in Q} d(x)$ and $d(0)~=~0$.
Also we assume that the function $\varphi(x)$ is strongly convex, i.e. $
\frac{\mu}{2}\|x~-~y\|^2~\leq~\varphi(y)~-~\varphi(x)~-~\la g(x) , y~-~x \ra
$ for all  $x,y \in Q, g(x) \in \partial \varphi(x)$. 
As a corollary, we have 
\begin{equation}
\varphi(x)-\varphi(x^*) \geq \frac{\mu}{2} \|x-x^*\|^2, \quad \forall x \in Q,
\label{eq:PhiStrConv}
\end{equation}
where $x^*$ is the solution of the problem \eqref{eq1:PrStateInit}.
We  also assume that $d(x)$ satisfies the following property. If $x_0$ is a random vector such that $\E_{x_0} \|x-x_0\|^2 \leq R_0^2$ for some fixed point $x$ and number $R_0$, then, for some $V >0$, 
\begin{equation}
\E_{x_0} d \left( \frac{x-x_0}{R_0} \right) \leq \frac{V^2}{2}.
\label{eq:dRestr}
\end{equation}

\begin{algorithm}[h!]
\caption{Stochastic Intermediate Gradient Method for Strongly Convex Problems}
\label{alg:SIGMA1}
\begin{algorithmic}[1]
\REQUIRE The function $d(x)$, point $u_0$, number $R_0$ such that $\|u_0-x^*\| \leq R_0$, number $p \in [1,2]$.
\ENSURE The point $u_{k+1}$.
\STATE Set $k$ = 0.
\STATE Calculate
\begin{equation}
N_k := \left\lceil \left(\frac{4\ee C_1LV^2}{\mu}\right)^{\frac{1}{p}}\right\rceil  .
\label{eq:SIGMA1N_0}
\end{equation}
\REPEAT
\STATE Calculate
\begin{align}
& m_k := \max \left\{1, \left\lceil \frac{16\ee^{k+2}  C_2^2 \sigma^2 V^2}{\mu^2 R_0^2 N_k } \right\rceil \right\},  \label{eq:SIGMA1m_k} \\
& R_k^2 := R_0^2 \ee^{-k} + \frac{2^p \ee C_3 \delta} {\mu(\ee-1)} \left(\frac{4\ee C_1LV^2}{\mu}\right)^{\frac{p-1}{p}} \left(1-\ee^{-k} \right). \label{eq:SIGMA1R_k}
\end{align}
\STATE Run Algorithm \ref{alg:SIGM} with $x_0=u_k$ and prox-function $d \left( \frac{x-u_k}{R_k} \right)$ for $N_k$ steps, using oracle $\tilde{G}^k_{\delta,L}(x) := \frac{1}{m_k}\sum_{i=1}^{m_k}{G_{\delta,L}(x,\xi^i)}$, where $\xi^i$, $i=1,...,m_k$ are i.i.d, on each step and sequences $\{\alpha_i\}_{i \geq 0}$, $\{\beta_i\}_{i \geq 0}$, $\{B_i\}_{i \geq 0}$ defined in Theorem \ref{Th04}. 
\STATE Set $u_{k+1}=y_{N_k}$, $k=k+1$.
\UNTIL{}
\end{algorithmic}
\end{algorithm}

\begin{Th}
After $k \geq 1$ outer iterations of Algorithm \ref{alg:SIGMA1}, we have
\begin{align}
&\E \varphi(u_{k}) - \varphi^* \leq \frac{\mu R_0^2}{2} \ee^{-k} + \frac{C_3 \ee 2^{p-1}}{\ee-1} \left(\frac{4\ee C_1LV^2}{\mu}\right)^{\frac{p-1}{p}} \delta, \label{eq:SIGMArate_phi} \\
& \E \|u_{k} - x^*\|^2 \leq R_0^2 \ee^{-k} + \frac{C_3 \ee 2^{p}}{\mu(\ee-1)} \left(\frac{4\ee C_1LV^2}{\mu}\right)^{\frac{p-1}{p}}  \delta. \label{eq:SIGMArate_x}
\end{align}
As a consequence, if we choose the error $\delta$ of the oracle satisfying
\begin{equation}
\delta \leq \frac{\e (\ee-1)}{2^pC_3\ee} \left(\frac{4\ee C_1LV^2}{\mu}\right)^{\frac{1-p}{p}},
\label{eq:SIGMADelta}
\end{equation}
then we need $N= \left\lceil \ln \left( \frac{\mu R_0^2}{ \e }\right) \right\rceil$ outer iterations and no more than
$$
\left(1+ \left(\frac{4\ee C_1LV^2}{\mu}\right)^{\frac{1}{p}} \right) \left( 1+ \ln \left( \frac{\mu R_0^2}{ \e }\right) \right) + \frac{16\ee^3C_2^2\sigma^2V^2}{\mu \e (\ee-1)}
$$
oracle calls to guarantee that $\E \varphi(u_{N}) - \varphi^* \leq \e$. 
\label{Th06}
\end{Th} 
To obtain complexity in terms of large deviations probability, we assume that the prox-function has quadratic growth with parameter $V^2$ with respect to the chosen norm, i.e.
\begin{equation}
d(x) \leq \frac{V^2}{2}\|x\|^2, \quad \forall x \in \R^n. 
\label{eq:dQuadrGrowth}
\end{equation}

Now we present a modification of Algorithm \ref{alg:SIGMA1} and a theorem with a bound for the probability of large deviations for the non-optimality gap of this algorithm. 
\begin{algorithm}[h!]
\begin{algorithmic}[1]
\REQUIRE The function $d(x)$, point $u_0$, number $R_0$ such that $\|u_0-x^*\| \leq R_0$, number $p \in [1,2]$, number $N \geq 1$ of outer iterations, confidence level $\Lambda$.
\ENSURE The point $u_{N}$.
\STATE Set $k$ = 0.
\STATE Calculate 
\begin{equation}
N_k := \left\lceil \left(\frac{6\ee C_1LV^2}{\mu}\right)^{\frac{1}{p}}\right\rceil  .
\label{eq:SIGMA2N_0}
\end{equation}
\STATE
\REPEAT
\STATE Calculate
\begin{align}
& m_k := \max \left\{1, \left\lceil \frac{36\ee^{k+2}  C_2^2 \sigma^2 V^2 \left(1+\ln \left(\frac{3N}{\Lambda} \right) \right)^2}{\mu^2 R_0^2 N_k } \right\rceil, \left\lceil \frac{144\ee^{k+2}  C_4^2 \sigma^2 \ln \left(\frac{3N}{\Lambda} \right) }{\mu^2 R_0^2 N_k } \right\rceil \right\},  \label{eq:SIGMA2m_k} \\
& R_k^2 := R_0^2 \ee^{-k} + \frac{2^p \ee C_3 \delta} {\mu(\ee-1)} \left(\frac{6\ee C_1LV^2}{\mu}\right)^{\frac{p-1}{p}} \left(1-\ee^{-k} \right), \label{eq:SIGMA2R_k} \\
& Q_k := \left\{ x \in Q : \|x-u_k\|^2 \leq R_k^2 \right\}. \label{eq:SIGMA2Q_k}
\end{align}
\STATE Run Algorithm \ref{alg:SIGM} applied to the problem $\min_{x \in Q_k} \varphi(x)$ with $x_0=u_k$ and prox-function $d \left( \frac{x-u_k}{R_k} \right)$ for $N_k$ steps using oracle $\tilde{G}^k_{\delta,L}(x) := \frac{1}{m_k}\sum_{i=1}^{m_k}{G_{\delta,L}(x,\xi^i)}$, where $\xi^i$, $i=1,...,m_k$ are i.i.d, on each step and sequences $\{\alpha_i\}_{i \geq 0}$, $\{\beta_i\}_{i \geq 0}$, $\{B_i\}_{i \geq 0}$ defined in Theorem \ref{Th04}. 
\STATE Set $u_{k+1}=y_{N_k}$, $k=k+1$.
 \UNTIL{$k=N-1$   }
\caption{Stochastic Intermediate Gradient Method for Strongly Convex Problems 2}
\label{alg:SIGMA2}
\end{algorithmic}
\end{algorithm}

\begin{Th}
After $N$ outer iterations of Algorithm \ref{alg:SIGMA2}, we have
\begin{align}
&\Prob \left\{ \varphi(u_{N}) - \varphi^* >  \frac{\mu R_0^2}{ 2 } \ee^{-N} + \frac{2^{p-1} \ee C_3 \delta} {(\ee-1)} \left(\frac{6\ee C_1LV^2}{\mu}\right)^{\frac{p-1}{p}} \delta \right\} \leq \Lambda. \label{eq:SIGMA2eL}
\end{align}
As a consequence, if we choose error of the oracle $\delta$ satisfying
\begin{equation}
\delta \leq \frac{\e (\ee-1)}{2^{p}C_3\ee} \left(\frac{6\ee C_1LV^2}{\mu}\right)^{\frac{1-p}{p}},
\label{eq:SIGMA2Delta}
\end{equation}
then we need no more than $N= \left\lceil \ln \left( \frac{\mu R_0^2}{ \e }\right) \right\rceil$ outer iterations and no more than
\begin{align}
&\left(1+ \left(\frac{6\ee C_1LV^2}{\mu}\right)^{\frac{1}{p}} \right) \left( 1+ \ln \left( \frac{\mu R_0^2}{ \e }\right) \right) + \notag \\
&+ \frac{36\ee^{3}  C_2^2 \sigma^2 V^2 }{\mu(\ee-1) \e } \left(1+\ln \left(\frac{3}{\Lambda}\left(1+\ln \left( \frac{\mu R_0^2}{ \e }\right) \right) \right) \right)^2 + \notag \\ 
& +\frac{144\ee^{3}  C_4^2 \sigma^2 }{\mu \e (\ee-1) } \ln \left(\frac{3}{\Lambda} \left(1+\ln \left( \frac{\mu R_0^2}{ \e }\right)\right) \right)  \label{eq:SIGMA2Compl}
\end{align}
oracle calls to guarantee that $\Prob \{ \varphi(u_{N}) - \varphi^* > \e \} \leq \Lambda$. 
\label{Th07}
\end{Th}

\subsection{Learning supervised pagerank with gradient-based and  gradient-free optimization methods.}
In this subsection we consider a parametric model for web-page ranking and learning the parameters of this model in a supervised setting. The results of this subsection are published in \cite{bogolubsky2016learning}.

\subsubsection{Model description}
\label{model}

Let $\Gamma=(V,E)$ be a directed graph. 
We suppose that for any $i\in V$ and any $i \rightarrow j\in E$, a vector of node's features $\mathbf{V}_i\in\mathbb{R}^{m_1}_+$ and a vector of edge's features $\mathbf{E}_{ ij}\in\mathbb{R}^{m_2}_+$ are given. Let $\vp_1 \in \R^{m_1}$, $\vp_2 \in \R^{m_2}$ be two vectors of parameters.
We denote $m=m_1+m_2$, $p=|V|$, $\vp=(\vp_1,\vp_2)^T$. 
Let us describe the random walk on the graph $\Gamma$.
A surfer starts a random walk from a random page $i\in U$ ($U$ is some subset in $V$ called {\it seed set}, $|U|=n$).
We assume that $\vp_1$ and node features are chosen in such way that $\sum_{l\in U} \langle\varphi_1,\mathbf{V}_{l}\rangle $ is non-zero. The initial probability of being at vertex $i \in V$ is called the {\it restart probability} and equals
\begin{equation}
 [\pi^0(\vp)]_i=\frac{\langle\varphi_1,\mathbf{V}_{i}\rangle}{\sum_{l\in  U}\langle\varphi_1,\mathbf{V}_{l}\rangle }, \quad i \in U
\label{restart}
\end{equation}
and $[\pi^0(\vp)]_i=0$ for $i \in V\setminus U$.
At each step, the surfer (with a current position $i\in V$) either chooses with probability $\alpha\in(0,1)$, 
which is called the {\it damping factor}, to go to any vertex from $V$ in accordance with the distribution $\pi^0(\vp)$ (makes a {\it restart}) or chooses to traverse an outgoing edge (makes a {\it transition}) with probability $1-\alpha$.
We assume that $\vp_2$ and edges features are chosen in such way that $\sum_{l: i \to  l}\langle\varphi_2,\mathbf{E}_{i l}\rangle$ is non-zero for all $i$ with non-zero outdegree.
For $i$ with non-zero outdegree, the probability
\begin{equation}
 [P(\vp)]_{i,j}=\frac{\langle\varphi_2,\mathbf{E}_{ij}\rangle }{\sum_{l: i \to  l}\langle\varphi_2,\mathbf{E}_{il}\rangle }
\label{transition}
\end{equation}
of traversing an edge $i\rightarrow j \in E$ is called the {\it transition probability}.
If an outdegree of $i$ equals $0$, then we set $[P(\varphi)]_{i,j}=[\pi^0(\varphi)]_j$ for all $j\in V$ (the surfer with current position $i$ makes a restart with probability $1$). Finally, by Equations \eqref{restart} and~\eqref{transition} the total probability of choosing vertex $j\in V$ conditioned by the surfer being at vertex $i$ equals $\alpha[\pi^0(\vp)]_j+(1-\alpha)[P(\vp)]_{i,j}$.
Denote by $\pi (\vp) \in \R^{p}$ the stationary distribution of the described Markov process.
It can be found as a solution of the system of equations
\begin{equation}
\pi = \alpha \pi^0(\vp)  + (1-\alpha) P^T(\varphi)  \pi
\label{eq:pi_Phi_P2}
\end{equation}
We learn the ranking algorithm, which orders the vertices $i$ by their probabilities $[\pi]_{i}$ in the stationary distribution $\pi$.

\subsubsection{Loss-minimization problem statement}
\label{optimal}

Let $Q$ be a set of queries and, for any $q\in Q$, a set of nodes $V_q$ which are relevant to $q$ be given. We are also provided with a ranking algorithm which assigns nodes ranking scores $[\pi_q]_i$, $i \in V_q$, $\pi_q=\pi_q(\varphi)$, as its output. For example, in web search, the score $[\pi_q]_i$ may repesent relevance of the page $i$ w.r.t. the query $q$. Our goal is to find the  parameter vector $\vp$ which minimizes the discrepancy of the ranking scores from the ground truth scoring defined by assessors. For each $q\in Q$, there is a set of nodes in $V_q$ manually judged and grouped by relevance labels $1,\ldots,\ell$. We denote $V^j_q$ the set of nodes annotated with label $\ell+1-j$ (i.e., $V_q^1$ is the set of all nodes with the highest relevance score).
Following the common approach, 
we consider the square
loss function and minimize
\begin{equation}
f(\varphi)= \frac{1}{|Q|} \sum_{q=1}^{|Q|} \|(A_q \pi_q (\vp))_{+} \|^2_2
\label{eq:f_phi_def_2}
\end{equation}
as a function of $\vp$ over some set of feasible values $\Phi$, where vector $x_+$ has components $[x_+]_i=\max\{x_i,0\}$, the matrices $A_q \in \R^{r_q \times p_q}, q \in Q$ represent assessor's view of the relevance of pages to the query $q$, $r_q$ equals $\sum_{1\leq j<l\leq \ell}|V_q^{j}||V_q^{l}|$. We denote $r=\max_{q \in Q} r_q$. 
By definition each row of matrix $A_q$ corresponds to some pair of pages $i_1 \in V_q^j$, $i_2 \in V_q^l$, where $j<l$, and the $i_1$-th element of this row is equal to $-1$, $i_2$-th element is equal to $1$, and all other elements are equal to $0$.

We consider the ranking algorithm based on scores~\eqref{eq:pi_Phi_P2} in Markov random walk on a graph $\Gamma_q=(V_q,E_q)$. We assume that feature vectors $\mathbf{V}_i^q$, $i\in V_q$,  $\mathbf{E}_{ij}^q$, $ i\rightarrow j\in E_q$, depend on $q$ as well. For example, vertices in $V_q$ may represent web pages which were visited by users after submitting a query $q$ and features may reflect different properties of query--page pair. For fixed $q\in Q$, we consider all the objects related to the graph $\Gamma_q$ introduced in the previous section: $U_q:=U$, $\pi_q^0:=\pi^0$, $P_q:=P$, $p_q:=p$, $n_q:=n$, $\pi_q:=\pi$. This allows ranking model to capture common (``static'') dependencies, which do vary between different queries. In this way, the ranking scores depend on query via the ``dynamic'' (query-dependent) features, but the parameters of the model $\alpha$ and $\vp$ are not query-dependent. 
We also denote $p=\max_{q \in Q} p_q$, $n=\max_{q \in Q} n_q$, $s=\max_{q \in Q} s_q$, where  $s_q=\max_{i \in V_q} |\{j: i \to j \in E_q\}|$. 
In order to guarantee that the probabilities in \eqref{restart} and \eqref{transition} are non-negative and that they do not blow up due to zero value of the denominator, we need appropriately choose the set $\Phi$ of possible values of parameters $\vp$. Thus we choose some $\hat{\vp}$ and $R >0 $ such that the set $\Phi$ defined as $\Phi = \{ \vp \in \R^m:  \|\vp-\hat{\vp}\|_2 \leq R \}$ lies in the set of vectors with positive components $\R^m_{++}$.
The loss-minimization problem which we solve is as follows
\begin{equation}
    \min_{ \vp \in \Phi} f(\vp), \Phi = \{ \vp \in \R^m:  \|\vp-\hat{\vp}\|_2 \leq R \}.
    \label{eq:prob_form}
\end{equation}

From \eqref{eq:pi_Phi_P2}, we obtain the following equation for $p_q \times m$ matrix $\frac{d \pi_q(\vp)}{d \vp^T}$ which is the derivative of stationary distribution $\pi_q(\vp)$ with respect to $\vp$
\begin{align}
& \frac{d \pi_q(\vp)}{d \vp^T} = \alpha \frac{d \pi^0_q(\vp)}{d \vp^T}  + (1-\alpha) \sum_{i=1}^{p_q} \frac{ d p_i(\vp)}{d \vp^T} [\pi_q(\vp)]_i + (1-\alpha) P_q^T(\vp) \frac{d \pi_q(\vp)}{d \vp^T},
\label{eq:pi_q_full_der}
\end{align}
where
$p_i(\vp)$ is the $i$-th column of the matrix $P_q^T(\vp)$. Then the gradient of the function $f(\vp)$ is easy to derive:
\begin{equation}
\nabla f(\vp) = \frac{2}{|Q|} \sum_{q=1}^{|Q|} \left(\frac{d \pi_q(\vp)}{d \vp^T} \right)^T A_q^T (A_q \pi_q (\vp))_{+}.
\label{eq:nf}
\end{equation}

The method~~\cite{nesterov2015finding} for approximation of $\pi_q(\vp)$ for any fixed $q\in Q$ constructs a sequence $\pi_k$ and the output $\tpi_q (\vp,N)$ (for some fixed non-negative integer $N$) 
by the following rule
\begin{equation}
\pi_0 = \pi_q^0(\vp), \quad \pi_{k+1} = P_q^T(\varphi) \pi_k, \quad \tpi_q (\vp,N) = \frac{\alpha}{1-(1-\alpha)^{N+1}} \sum_{k=0}^{N} {(1-\alpha)^k \pi_k}.
\label{eq:pi_k+1_def}
\end{equation}
\begin{Lm}
Assume that for some $\delta_1>0$ Method \eqref{eq:pi_k+1_def} with
$
N = \left\lceil \frac{1}{\alpha} \ln \frac{8r}{\delta_1} \right\rceil - 1
$
is used to calculate the vector $\tpi_q (\vp,N)$ for every $q\in Q$. Then
\begin{equation}
\tf(\varphi,\delta_1) = \frac{1}{|Q|} \sum_{q=1}^{|Q|} \|(A_q \tpi_q (\vp,N))_{+} \|^2_2
\label{eq:tfN1_def}
\end{equation}
satisfies
\begin{equation}
|\tf(\varphi,\delta_1) - f(\vp)|\leq \delta_1.
\label{eq:tf_error}
\end{equation}
Moreover, the calculation of $\tf(\varphi,\delta_1)$ requires not more than
$
|Q|(3 mps + 3 psN+6r)
$
a.o.
and not more than $3ps$ memory items.
\label{Lm:f_compl}
\end{Lm}

Our generalization of the method~\cite{nesterov2015finding} for calculation of $\frac{d \pi_q(\vp)}{d \vp^T}$ for any $q\in Q$  is the following.
Choose some non-negative integer $N_1$ and calculate $\tpi_q(\vp,N_1)$ using \eqref{eq:pi_k+1_def}. Calculate a sequence $\Pi_k$
\begin{equation}
\Pi_0 = \alpha \frac{d \pi^0_q(\vp)}{d \vp^T}  + (1-\alpha) \sum_{i=1}^{p_q} \frac{ d p_i(\vp)}{d \vp^T} [\tpi_q(\vp,N_1)]_i, \quad \Pi_{k+1} = P_q^T(\vp) \Pi_k.
\label{eq:Pi_0_def}
\end{equation}
The output is (for some fixed non-negative integer $N_2$)
\begin{equation}
\tPi_q(\vp,N_2) = \frac{1}{1-(1-\alpha)^{N_2+1}} \sum_{k=0}^{N_2} (1-\alpha)^k \Pi_k.
\label{eq:tPi_def}
\end{equation}

In what follows, we use the following norm on the space of matrices $A \in \R^{n_1\times n_2}$: $\|A\|_1 = \max_{j = 1,...,n_2} \sum_{i=1}^{n_1} |a_{ij}|$.

\begin{Lm}
Let $\beta_1$ be a number (explicitly computable, see \cite{bogolubsky2016learning}) such that for all $\vp \in \Phi$
\begin{equation}
\alpha \left\|\frac{d \pi^0_q(\vp)}{d \vp^T} \right\|_1 + (1-\alpha) \sum_{i=1}^{p_q} \left\|\frac{ d p_i(\vp)}{d \vp^T} \right\|_1 \leq  \beta_1.
\label{eq:Pi_0_est}
\end{equation}
Assume that Method \eqref{eq:pi_k+1_def} with
$
N_1= \left\lceil \frac{1}{\alpha} \ln \frac{24\beta_1 r}{\alpha \delta_2} \right\rceil - 1
$
is used for every $q \in Q$ to calculate the vector $\tpi_q(\vp,N_1)$ and Method \eqref{eq:Pi_0_def}, \eqref{eq:tPi_def}
with
$
N_2=\left\lceil \frac{1}{\alpha} \ln \frac{8\beta_1 r}{\alpha \delta_2} \right\rceil - 1
$
is used for every $q \in Q$ to calculate the matrix $\tPi_q(\vp,N_2)$ \eqref{eq:tPi_def}.
Then the vector
\begin{equation}
\tg(\vp,\delta_2) = \frac{2}{|Q|} \sum_{q=1}^{|Q|} \left(\tPi_q(\vp,N_2) \right)^T A_q^T (A_q \tpi_q(\vp,N_1)  )_{+}
\label{eq:tnfN2_def}
\end{equation}
satisfies
\begin{equation}
\left\| \tg(\vp,\delta_2) - \nabla f(\vp)\right\|_{\infty}  \leq \delta_2.
\label{eq:tnf_error}
\end{equation}
Moreover the calculation of $\tg(\vp,\delta_2)$ requires no more than
$
|Q| (10 mps + 3 psN_1+ 3mpsN_2  + 7r)
$
a.o. and not more than $4ps+4mp+r$ memory items.
\label{Lm:nf_compl}
\end{Lm}
As we see, there is an inexact oracle available for the considered supervised learning problem. Thus, in the next subsections, we consider a general problem with intexact oracle and solve it by zero-order and first-order methods.

\subsubsection{Solving the learning problem by zero-order method}
First, we consider a general zero-order method with inexact function evaluations and then we apply it to solve the learning problem. Let $\Ec$ be an $m$-dimensional vector space. First, we consider a general function $f(\cdot): \Ec \to \R$ and denote its argument by $x$ or $y$ to avoid confusion with the above text. We denote the value of linear function $g \in \Ec^*$ at $x\in \Ec$ by $\la g, x \ra$. We choose some norm $\|\cdot\|$ in $\Ec$ and say that $f \in C^{1,1}_{L} (\|\cdot\|)$ iff
\begin{equation}
     |f(x)-f(y) - \la \nabla f(y) ,x-y \ra| \leq \frac{L}{2}\|x-y\|^2, \quad \forall x,y \in \Ec.
     \label{eq:fLipSm}
\end{equation}
The problem of our interest is to find $\min_{x\in X} f(x)$, where $f \in C^{1,1}_{L}(\|\cdot\|)$, $X$ is a closed convex set and there exists a number $D \in (0,+\infty)$ such that ${\rm diam} X := \max_{x,y \in X}\|x-y\| \leq D$. Also we assume that the inexact zero-order oracle for $f(x)$ returns a value $\tf(x,\delta)=f(x) + \tilde{\delta}(x)$, where $\tilde{\delta}(x)$ is the error satisfying for some $\delta >0$ (which is known) $|\tilde{\delta}(x)| \leq \delta$ for all $x \in X$.
Let $x^* \in \arg \min_{x\in X} f(x)$. Denote $f^* = \min_{x\in X} f(x)$.

Unlike \cite{nesterov2017random}, we define the biased gradient-free oracle $g_{\tau}(x,\delta) = \frac{m}{\tau}(\tf(x+\tau \xi,\delta) -\tf(x,\delta)) \xi$,
where $\xi$ is a random vector uniformly distributed over the unit sphere $\S =\{ t \in \R^m  : \|t\|_2 = 1\}$, $\tau$ is a smoothing parameter.

Algorithm~\ref{alg:GFPGM} below is the variation of the projected gradient descent method. Here $\Pi_{X}(x)$ denotes the Euclidean projection of a point $x$ onto the set $X$.
\begin{algorithm}[h!]
        \caption{Gradient-type method}
        \label{alg:GFPGM}
\begin{algorithmic}[1]
   \STATE {\bfseries Input:} Point $x_0 \in X$, stepsize $h>0$, number of steps $M$.
   \STATE Set $k=0$.
   \REPEAT
   \STATE Generate $\xi_k$ and calculate corresponding $g_{\tau}(x_k,\delta)$.
   \STATE Calculate $x_{k+1} = \Pi_{X}(x_k- h g_{\tau}(x_k,\delta))$.
         \STATE Set $k=k+1$.
   \UNTIL{$k>M$}
         \STATE {\bfseries Output:} The point $y_M = \arg \min_x \{ f(x): x \in \{ x_0, \dots, x_M\}\}$.
\end{algorithmic}
\end{algorithm}

Next theorem gives the convergence rate of Algorithm \ref{alg:GFPGM}.
Denote by $\Xi_k=(\xi_0, \dots, \xi_k)$ the history of realizations of the vector $\xi$ generated on each iteration of the algorithm.

\begin{Th}
Let $f \in C^{1,1}_{L} (\|\cdot\|_2)$ and convex. Assume that  $x^*\in {\rm int} X$, and the sequence $x_k$ is generated by Algorithm \ref{alg:GFPGM} with $h=\frac{1}{8mL}$.
Then for any $M \geq 0$, we have
\begin{align}
& \E_{\Xi_{M-1}} f(y_M) - f^*\leq \frac{8mL D^2}{M+1} + \frac{\tau^2 L (m+8)}{8} + \frac{\delta m D}{4 \tau }  + \frac{\delta^2 m}{L \tau^2}.
\label{eq:rtSmth}
\end{align}
\label{th_1}
\end{Th}

\begin{algorithm}[h!]
        \caption{Gradient-free method for Problem \eqref{eq:prob_form}}
        \label{alg:RGFGM_learn}
\begin{algorithmic}[1]
   \STATE {\bfseries Input:} Point $\varphi_0 \in \Phi$, $L$ -- Lipschitz constant for the function $f(\varphi)$ on $\Phi$,    accuracy $\varepsilon >0$.
         \STATE Define $M=\left\lceil128 m \frac{LR^2}{\varepsilon}\right\rceil$, $\delta=  \frac{\varepsilon^{\frac32}\sqrt{2}}{16mR\sqrt{L(m+8)}} $, $\tau = \sqrt{\frac{ 2\varepsilon}{L (m+8)}}$.
         \STATE Set $k=0$.
   \REPEAT
    \STATE Generate random vector $\xi_k$ uniformly distributed over a unit Euclidean sphere $\S$ in $R^m$.
    \STATE Calculate $\tf(\varphi_k+ \tau \xi_k,\delta)$, $\tf(\varphi_k,\delta)$ using Lemma \ref{Lm:f_compl} with $\delta_1=\delta$.
    \STATE Calculate  $g_{\tau}(\vp_k,\delta) = \frac{m}{\tau}(\tf(\varphi_k+ \tau \xi_k,\delta)-\tf(\varphi_k,\delta))\xi_k$.
         \STATE Calculate $\vp_{k+1} = \Pi_{\Phi}\left(\vp_k- \frac{1}{8 m L} g_{\tau}(\vp_k,\delta)\right)$.
         \STATE Set $k=k+1$.
   \UNTIL{$k>M$}
         \STATE {\bfseries Output:} The point $\hat{\vp}_M = \arg \min_{\vp} \{ f(\vp): \vp \in \{ \vp_0, \dots, \vp_M\}\}$.
\end{algorithmic}
\end{algorithm}

Next, we apply the above method to solve the learning problem \eqref{eq:prob_form}. The resulting algorithm is listed as Algorithm \ref{alg:RGFGM_learn}.
The most computationally hard on each iteration of the main cycle of this method are calculations of $\tf(\varphi_k+ \tau \xi_k,\delta)$, $\tf(\varphi_k,\delta)$. Using Lemma \ref{Lm:f_compl}, we obtain that each iteration of Algorithm \ref{alg:RGFGM_learn} needs no more than
$$
2|Q|\left(3 mps + \frac{3 ps}{\alpha}\ln \frac{128mrR\sqrt{L(m+8)}}{\varepsilon^{3/2}\sqrt{2}} + 6r\right)
$$
a.o.  So, we obtain the following result, which gives the complexity of Algorithm \ref{alg:RGFGM_learn}.
\begin{Th}
Assume that the set $\Phi$ in \eqref{eq:prob_form} is chosen in a way such that $f(\varphi)$ is convex on $\Phi$ and some $\vp^* \in \arg \min_{\vp \in \Phi}f(\vp)$ belongs also to $ {\rm int} \Phi$. Then the mean total number of arithmetic operations of the Algorithm \ref{alg:RGFGM_learn} for the accuracy $\e$ (i.e. for the inequality $\E_{\Xi_{M-1}} f(\hat{\vp}_M) - f(\vp^*) \leq \e$ to hold) is no more than
$$
768 mps|Q| \frac{LR^2}{\varepsilon} \left( m  + \frac{1}{\alpha}  \ln \frac{128mrR\sqrt{L(m+8)}}{\varepsilon^{3/2}\sqrt{2}}+6r\right).
$$
\label{Th:RGFGM_learn_compl}
\end{Th}

\subsubsection{Solving the learning problem by first-order method}
First we consider a general first-order method with inexact function values and inexact gradient, and then we apply it to solve the learning problem.
We generalize the approach in~\cite{ghadimi2016mini-batch} for constrained non-convex optimization problems. Our main contribution consists in developing this framework for an inexact first-order oracle and unknown "Lipschitz constant" of this oracle.

Let $\Ec$ be a finite-dimensional real vector space and $\Ec^*$ be its dual. We denote the value of linear function $g \in \Ec^*$ at $x\in \Ec$ by $\la g, x \ra$. Let $\|\cdot\|$ be some norm on $\Ec$, $\|\cdot\|_*$ be its dual.
Our problem of interest in this subsection is a {\it composite optimization} problem of the form
\begin{equation}
\min_{x \in X} \{ \psi(x) := f(x) + h(x)\},
\label{eq:PrStateInit}
\end{equation}
where $X \subset \Ec$ is a closed convex set, $h(x)$ is a simple convex function, e.g. $\|x\|_1$. We assume that $f(x)$ is a general function endowed with an inexact first-order oracle in the following sense. There exists a number $L \in (0,+\infty)$ such that for any $\delta \geq 0$ and any $x\in X$ one can calculate $\tf(x,\delta) \in \R$ and $\tg(x,\delta) \in \Ec^*$ satisfying
\begin{equation}
     |f(y)-(\tf(x,\delta)  - \la\tg(x,\delta) ,y-x \ra)| \leq \frac{L}{2}\|x-y\|^2 + \delta.
     \label{eq:dL_or_def}
\end{equation}
for all $y \in X$. The constant $L$ can be considered as "Lipschitz constant" because for the exact first-order oracle for a function $f \in C_L^{1,1} (\|\cdot\|)$  \eqref{eq:dL_or_def} holds with $\delta = 0$. This is a generalization of the concept of $(\delta,L)$-oracle considered in \cite{devolder2013exactness} for convex problems.

We choose a {\it prox-function} $d(x)$ which is continuously differentiable and $1$-strongly convex on $X$ with respect to $\|\cdot\|$.
This means that for any $x,y \in X$ $d(y)-d(x) -\la \nabla d(x) ,y-x \ra \geq \frac12\|y-x\|^2$.
We define also the corresponding {\it Bregman distance} $V (x, z) = d(x) - d(z) - \la \nabla d(z), x - z \ra$.

%

\begin{algorithm}[h!]
        \caption{Adaptive projected gradient algorithm}
        \label{alg:gen_PG_LA_2}
\begin{algorithmic}[1]
   \STATE {\bfseries Input:} Point $x_0 \in X$, number $L_0 >0$.
   \STATE Set $k=0$, $z=+\infty$.
   \REPEAT
			\STATE Set $M_k=L_k$, ${\rm flag}=0$.
			\REPEAT
				\STATE Set $\delta = \frac{\e}{16M_k}$. Calculate $\tf(x_k,\delta)$ and $\tg(x_k,\delta)$.
				\STATE $
				w_k
				= \arg \min_{x \in Q} \left\{\la \tg(x_k,\delta),x \ra + M_kV(x,x_k) +h(x) 		\right\}$
				\STATE If the inequality
				\begin{align*}
				&\tf(w_k,\delta) \leq \tf(x_k,\delta) + \la \tg(x_k,\delta) ,w_k - x_k \ra  +\frac{M_k}{2}\|w_k - x_k\|^2 +	\frac{\e}{8M_k}
				\end{align*}
				holds, set ${\rm flag}=1$. Otherwise set $M_k=2M_k$.
			\UNTIL{${\rm flag}=1$}
			\STATE Set $x_{k+1} = w_k$, $L_{k+1}=\frac{M_k}{2}$.
			\STATE If $\left\|M_k(x_k-x_{k+1})\right\| < z$, set $z=\left\|M_k(x_k-x_{k+1})\right\|$, $K=k$.
      \STATE Set $k=k+1$.
   \UNTIL{$z \leq \e$}
         \STATE {\bfseries Output:} The point $x_{K+1}$.
\end{algorithmic}
\end{algorithm}

\begin{Th}
Assume that $f(x)$ is endowed with the inexact first-order oracle in the sense of \eqref{eq:dL_or_def} and that there exists a number $\psi^* > -\infty$ such that $\psi(x) \geq \psi^*$ for all $x \in X$. Then after $M$ iterations of Algorithm \eqref{alg:gen_PG_LA_2} it holds that
\begin{equation}
\left\|M_K(x_K-x_{K+1})\right\|^2 \leq \frac{4L(\psi(x_0)-\psi^*)}{M+1} + \frac{\e}{2}.
\label{eq:pg_la_2_rate}
\end{equation}
Moreover, the total number of inner steps is no more than $M+\log_2\frac{2L}{L_0}$.
\label{Th:pg_la_2_rate}
\end{Th}

Next we apply the general method to the learning problem.
Since $\Phi$ is the Euclidean ball, it is natural to set $\Ec = R^m$ and $\|\cdot\|= \|\cdot\|_2$, choose the prox-function $d(\vp)=\frac12\|\vp\|_2^2$. Then the Bregman distance is $V(\vp,\omega)=\frac{1}{2}\|\vp-\omega\|_2^2$.
Algorithm \ref{alg:NCGM_learn} is a formal record of the algorithm. 

\begin{algorithm}[h!]
        \caption{Adaptive gradient method for Problem~\eqref{eq:prob_form}}
        \label{alg:NCGM_learn}
\begin{algorithmic}[1]
   \STATE {\bfseries Input:} Point $\vp_0 \in \Phi$, number $L_0 >0$, accuracy $\varepsilon >0$.
    \STATE Set $k=0$, $z=+ \infty$.
		\REPEAT
			\STATE Set $M_k=L_k$, ${\rm flag}=0$.
			\REPEAT
				\STATE Set $\delta_1 = \frac{\e}{32M_k}$, $\delta_2 = \frac{\e}{64M_kR\sqrt{m}}$.
				\STATE Calculate $\tf(\vp_k, \delta_1)$ using Lemma \ref{Lm:f_compl} and $\tg(\vp_k, \delta_2)$ using Lemma \ref{Lm:nf_compl}.
				\STATE Find
				\begin{equation}
				\omega_k= \arg \min_{\vp \in \Phi} \left\{\la \tg(\vp_k, \delta_2),\vp \ra + \frac{M_k}{2}\|\vp-\vp_k\|_2^2. \right\}
			  \notag
				\end{equation}
				\STATE Calculate $\tf(\omega_k, \delta_1)$ using Lemma \ref{Lm:f_compl}.
				\STATE If the inequality
				\begin{align}
				&\tf(\omega_k, \delta_1) \leq 	\tf(\vp_k, \delta_1) + \la \tg(\vp_k, \delta_2) ,\omega_k - \vp_k \ra  +\frac{M_k}{2}\|\omega_k - \vp_k\|_2^2 		+		\frac{\e}{8M_k}
				\notag
				\end{align}
				holds, set ${\rm flag}=1$. Otherwise set $M_k=2M_k$.
			\UNTIL{${\rm flag}=1$}
			\STATE Set $\vp_{k+1} = \omega_k$, $L_{k+1}=\frac{M_k}{2}$, .
			\STATE If $\left\|M_k(\vp_k-\vp_{k+1})\right\|_2 < z$, set $z=\left\|M_k(\vp_k-\vp_{k+1})\right\|_2$, $K=k$.
      \STATE Set $k=k+1$.
   \UNTIL{$z\leq \e$}				
         \STATE {\bfseries Output:} The point $\vp_{K+1}$.
\end{algorithmic}
\end{algorithm}

The most computationally consuming operations of the inner cycle of Algorithm \ref{alg:NCGM_learn} are calculations of $\tf(\vp_k, \delta_1)$, $\tf(\omega_k, \delta_1)$ and $\tg(\vp_k, \delta_2)$. Using Lemma \ref{Lm:f_compl} and Lemma \ref{Lm:nf_compl}, we obtain that each inner iteration of Algorithm \ref{alg:NCGM_learn} needs no more than
$$
7r|Q| + \frac{6mps|Q|}{\alpha} \ln \frac{1024\beta_1 r R L \sqrt{m}}{ \alpha \e }
$$
a.o.
Using Theorem \ref{Th:pg_la_2_rate}, we obtain the following result, which gives the complexity of Algorithm \ref{alg:NCGM_learn}.
\begin{Th}
The total number of arithmetic operations in Algorithm \ref{alg:NCGM_learn} for the accuracy $\e$ (i.e. for the inequality $\left\|M_K(\vp_K-\vp_{K+1})\right\|_2^2 \leq \e$ to hold) is no more than
\begin{align}
    &\left(\frac{8L(f(\vp_0)-f^*)}{\e}+\log_2\frac{2L}{L_0}\right) \cdot  \left(7r|Q| + \frac{6mps|Q|}{\alpha} \ln \frac{1024\beta_1 r R L \sqrt{m}}{ \alpha \e } \right). \notag
\end{align}
\label{Th:NCGM_learn_compl}
\end{Th}

\subsection{An accelerated directional derivative method for smooth stochastic convex optimization.}
In this section we consider directional derivatives methods with inexact oracle for stochastic convex optimization. The results of this subsection are published in \cite{dvurechensky2020accelerated}.
Motivated by potential presence of non-stochastic noise in an optimization problem,
we assume that the noise in the directional derivative consists of two parts. Similar to stochastic optimization problems, the first part is of a stochastic nature. On the opposite, the second part is an additive noise of an unknown nature, but bounded in the absolute value. 
More precisely, we consider the following optimization problem
\begin{equation}
\label{eq:PrSt}
   \min_{x\in\R^n} \left\{ f(x) := \EE_{\xi}[F(x,\xi)] = \int_{\mathcal{X}}F(x,\xi)dP(x) \right\}, 
\end{equation}
where $\xi$ is a random vector with probability distribution $P(\xi)$, $\xi \in \mathcal{X}$, and for $P$-almost every $\xi \in \mathcal{X}$, the function $F(x,\xi)$ is closed  and  convex. Moreover, we assume that, for $P$ almost every $\xi$, the function $F(x,\xi)$ has gradient $g(x,\xi)$, which is $L(\xi)$-Lipschitz continuous with respect to the Euclidean norm and there exists $L_2\geqslant 0$ such that $\sqrt{\EE_{\xi} L(\xi)^2 } \leqslant L_2 < +\infty$. Under this assumptions, $\EE_{\xi}g(x,\xi) = \nabla f(x)$ and $f$ has $L_2$-Lipschitz continuous gradient with respect to the Euclidean norm. Also we assume that
\begin{equation}
\label{stoch_assumption_on_variance}
    \EE_{\xi}[\|g(x,\xi) - \nabla f(x)\|_2^2] \leqslant \sigma^2,
\end{equation}
where $\|\cdot\|_2$ is the Euclidean norm.

Finally, we assume that an optimization procedure, given a point $x \in \R^{n}$, direction $e \in S_2(1)$ and $\xi$ independently drawn from $P$, can obtain a noisy stochastic approximation $\tf'(x,\xi,e)$ for the directional derivative $\la g(x,\xi),e \ra$: 
\begin{align}
\tf'(x,\xi,e) &= \la g(x,\xi),e \ra + \zeta(x,\xi,e) + \eta(x,\xi,e), \notag \\ \EE_{\xi}(\zeta(x,\xi,e))^2 &\leqslant \Delta_{\zeta}, \;\forall x \in \R^n, \forall e \in S_2(1),  \notag \\ |\eta(x,\xi,e)| & \leqslant  \Delta_{\eta}, \;\forall x \in \R^n, \forall e \in S_2(1), \; \text{a.s. in } \xi,
\label{eq:tf_def}
\end{align}
where $S_2(1)$ is the Euclidean sphere of radius one with the center at the point zero and the values $\Delta_{\zeta}$, $\Delta_{\eta}$ are controlled and can be made as small as it is desired. Note that we use the smoothness of $F(\cdot,\xi)$ to write the directional derivative as $\la g(x,\xi),e \ra$, but we \textit{do not assume} that the whole stochastic gradient $g(x,\xi)$ is available.
We choose a \textit{prox-function} $d(x)$ which is continuous, convex on $\R^n$ and is $1$-strongly convex on $\R^n$ with respect to $\|\cdot\|_p$, i.e., for any $x, y \in \R^n$ $d(y)-d(x) -\la \nabla d(x) ,y-x \ra \geq \frac12\|y-x\|_p^2$.
Without loss of generality, we assume that $\min\limits_{x\in \R^n} d(x) = 0$.
 We define also the corresponding \textit{Bregman divergence} $V[z] (x) = d(x) - d(z) - \la \nabla d(z), x - z \ra$, $x, z \in \R^n$. 
For the case $p=1$, we choose the following prox-function \cite{ben-tal2015lectures}
\begin{equation}
\label{eq:dp1}
d(x) = \frac{{\rm e}n^{(\kappa-1)(2-\kappa)/\kappa}\ln n}{2} \|x\|_\kappa^2, \quad \kappa=1 + \frac{1}{\ln n}
\end{equation}
and, for the case $p=2$, we choose the prox-function to be the squared Euclidean norm $d(x) = \frac{1}{2}\|x\|_2^2$.

Based on the noisy stochastic observations \eqref{eq:tf_def} of the directional derivative, we form the following stochastic approximation of $\nabla f(x)$
\begin{equation}
\label{eq:MiniBatcStocGrad}
\tnmf(x) = \frac{1}{m}\sum\limits_{i=1}^m\tf'(x,\xi_i,e)e,
\end{equation}
where $e \in RS_2(1)$, $\xi_i$, $i=1,...,m$ are independent realizations of $\xi$, $m$ is the \textit{batch size}.

\subsubsection{Algorithms and main results for convex problems}
Our Accelerated Randomized Directional Derivative (ARDD) method is listed as Algorithm~\ref{Alg:ARDS}.
\begin{algorithm}
	\caption{Accelerated Randomized Directional Derivative (ARDD) method}
	\label{Alg:ARDS}
	\begin{algorithmic}[1]
		\REQUIRE $x_0$~---starting point; $N \geqslant 1$~--- number of iterations; $m \geqslant 1$~--- batch size.
		\ENSURE point $y_N$.
		\STATE $y_0 \leftarrow x_0, \, z_0 \leftarrow x_0$.
		\FOR{$k=0,\, \dots, \, N-1$.}
		\STATE $\alpha_{k+1} \leftarrow \frac{k+2}{96n^2\rho_n L_2}, \, \tau_{k} \leftarrow \frac{1}{48\alpha_{k+1}n^2\rho_n L_2} = \frac{2}{k+2}$.
		\STATE Generate $e_{k+1} \in RS_2(1)$ independently from previous iterations and $\xi_i$, $i=1,...,m$ -- independent realizations of $\xi$. 
		\STATE 
		$
		\tnmf(x_{k+1})= \frac{1}{m}\sum\limits_{i=1}^m\tf'(x_{k+1},\xi_i,e)e.
		$
		\STATE $x_{k+1} \leftarrow \tau_kz_k + (1-\tau_k)y_k $.
		\STATE $y_{k+1} \leftarrow x_{k+1}-\frac{1}{2L_2}\tnmf(x_{k+1})$.
		\STATE $z_{k+1} \leftarrow \argmin\limits_{z\in\R^n} \left\{ {\alpha_{k+1} n \left\langle \tnmf(x_{k+1}), \, z-z_k \right\rangle +V[z_{k}] \left( z \right)}\right\}$.
		\ENDFOR
		\RETURN $y_N$
	\end{algorithmic}
\end{algorithm}
%
\begin{Th}
\label{Th:ARDFDSConv}
    Let ARDD method be applied to solve problem \eqref{eq:PrSt}. Then
    \begin{equation}\label{eq:ARDFDSConv}
        \begin{array}{rl}
        \EE[f(y_N)] - f(x^*)
        \leqslant \frac{384\Theta_p n^2\rho_nL_2}{N^2} 
        + \frac{4N}{nL_2}\cdot\frac{\sigma^2}{m} + \frac{61N}{24L_2}\Delta_\zeta + \frac{122N}{3L_2}\Delta_\eta^2 \\
        + \frac{12\sqrt{2n\Theta_p}}{N^2}\left(\frac{\sqrt{\Delta_\zeta}}{2}+ 2\Delta_\eta\right)
        + \frac{N^2}{12n\rho_nL_2} \left(\frac{\sqrt{\Delta_\zeta}}{2} 
        + 2\Delta_\eta\right)^2,
    \end{array}
    \end{equation}
    where $\Theta_p = V[z_0](x^*)$ is defined by the chosen proximal setup and  $\EE[\cdot] = \EE_{e_1,\ldots,e_N,\xi_{1,1},\ldots,\xi_{N,m}}[\cdot]$.
\end{Th}
Before we proceed to the non-accelerated method, we give the appropriate choice of the ARDD method parameters $N$, $m$, and accuracy of the directional derivative evaluation $\Delta_\zeta$, $\Delta_\eta$. These values are chosen such that the r.h.s. of \eqref{eq:ARDFDSConv} is smaller than $\e$. For simplicity we omit numerical constants and summarize the obtained values of the algorithm parameters in Table 1 below. The last row represents the total number $Nm$ of oracle calls.  
{\footnotesize
\renewcommand{\arraystretch}{2}
\begin{table}[h]
    \hspace{-3em}\begin{tabular}{|c|c|c|}
        \hline
         & $p=1$ & $p=2$ \\
        \hline
        $N$ & $\sqrt{\frac{ n\ln nL_2 \Theta_1}{\varepsilon}}$ & $\sqrt{\frac{ n^{2}L_2\Theta_2}{\varepsilon}}$ \\
        \hline
        $m$ & $\max\left\{1,\sqrt{\frac{\ln n}{n}}\cdot\frac{\sigma^2}{\varepsilon^{3/2}}\cdot\sqrt{\frac{\Theta_1}{L_2}}\right\}$ & $\max\left\{1,\frac{\sigma^2}{\varepsilon^{3/2}}\cdot\sqrt{\frac{\Theta_2}{L_2}}\right\}$\\
        \hline
        $\Delta_\zeta$ & $\min\left\{n(\ln n)^2L_2^2\Theta_1,\, \frac{\varepsilon^2}{n\Theta_1},\, \frac{\varepsilon^{\frac{3}{2}}}{\sqrt{n\ln n}}\cdot\sqrt{\frac{L_2}{\Theta_1}}\right\}$ & $\min\left\{n^3L_2^2\Theta_2,\, \frac{\varepsilon^2}{n\Theta_2},\, \frac{\varepsilon^{\frac{3}{2}}}{n}\cdot\sqrt{\frac{L_2}{\Theta_2}}\right\}$\\
        \hline
        $\Delta_\eta$ & $\min\left\{\sqrt{n}\ln nL_2\sqrt{\Theta_1},\, \frac{\varepsilon}{\sqrt{n\Theta_1}},\, \frac{\varepsilon^{\frac{3}{4}}}{\sqrt[4]{n\ln n}}\cdot\sqrt[4]{\frac{L_2}{\Theta_1}}\right\}$ & $\min\left\{n^{\frac{3}{2}}L_2\sqrt{\Theta_2},\, \frac{\varepsilon}{\sqrt{n\Theta_2}},\, \frac{\varepsilon^{\frac{3}{4}}}{\sqrt{n}}\cdot\sqrt[4]{\frac{L_2}{\Theta_2}}\right\}$\\
        \hline
        Calls & $\max\left\{\sqrt{\frac{ n\ln nL_2\Theta_1}{\varepsilon}}, \frac{\sigma^2\Theta_1 \ln n}{\varepsilon^2}\right\}$ & $\max\left\{\sqrt{\frac{ n^{2}L_2\Theta_2}{\varepsilon}}, \frac{\sigma^2\Theta_2 n}{\varepsilon^2}\right\})$ \\
        \hline
    \end{tabular}
    \caption{Algorithm~\ref{Alg:ARDS} parameters for the cases $p=1$ and $p=2$.}
    \label{tab:ARDD}
\end{table}
}

Our Randomized Directional Derivative (RDD) method is listed as Algorithm \ref{Alg:RDFDS}.
\begin{algorithm}
	\caption{Randomized Directional Derivative (RDD) method}
	\label{Alg:RDFDS}
	\begin{algorithmic}[1]
		\REQUIRE $x_0$~---starting point; $N \geqslant 1$~--- number of iterations; $m \geqslant 1$~--- batch size.
		\ENSURE point $\bar{x}_N$.
		\FOR{$k=0,\, \dots, \, N-1$.}
		\STATE $\alpha \leftarrow \frac{1}{48n\rho_n L_2}$.
		\STATE Generate $e_{k+1} \in RS_2\left( 1 \right)$ independently from previous iterations and $\xi_i$, $i=1,...,m$ -- independent realizations of $\xi$. 
		\STATE $
		\tnmf(x_{k})= \frac{1}{m}\sum\limits_{i=1}^m\tf'(x_{k},\xi_i,e)e.
		$
		\STATE $x_{k+1} \leftarrow \argmin\limits_{x\in\R^n} \left\{ {\alpha n \left\langle \tnmf(x_{k}), \, x-x_k \right\rangle +V[x_{k}] \left( x \right)}\right\}$.
		\ENDFOR
		\RETURN $\bar{x}_N \leftarrow \frac{1}{N}\sum\limits_{k=0}^{N-1}x_k$
	\end{algorithmic}
\end{algorithm}

\begin{Th}\label{theorem_convergence_mini_gr_free_non_acc}
    Let RDD method 
    be applied to solve problem \eqref{eq:PrSt}. Then
    \begin{align}
        &\EE[f(\bar{x}_N)] - f(x_*) \leqslant \frac{384n\rho_nL_2\Theta_p}{N} + \frac{2}{L_2}\frac{\sigma^2}{m}+ \frac{n}{12L_2}\Delta_\zeta + \frac{4n}{3L_2}\Delta_\eta^2 \notag \\
        &+ \frac{8\sqrt{2n\Theta_p}}{N}\left(\frac{\sqrt{\Delta_\zeta}}{2} + 2\Delta_\eta\right) + \frac{N}{3L_2\rho_n}\left(\frac{\sqrt{\Delta_\zeta}}{2} + 2\Delta_\eta\right)^2, \label{theo_main_result_mini_gr_free}
    \end{align}
     where $\Theta_p = V[z_0](x^*)$ is defined by the chosen proximal setup and  $\EE[\cdot] = \EE_{e_1,\ldots,e_N,\xi_{1,1},\ldots,\xi_{N,m}}[\cdot]$.
\end{Th}
Before we proceed, we give the appropriate choice of the RDD method parameters $N$, $m$, and accuracy of the directional derivative evaluation $\Delta_\zeta$, $\Delta_\eta$. These values are chosen such that the r.h.s. of \eqref{theo_main_result_mini_gr_free} is smaller than $\e$. For simplicity we omit numerical constants and summarize the obtained values of the algorithm parameters in Table 2 below. The last row represents the total number $Nm$ of oracle calls, that is, the number of directional derivative evaluations. 
{
\renewcommand{\arraystretch}{2}
\begin{table}[h]
    \centering
    \begin{tabular}{|c|c|c|}
        \hline
         & $p=1$ & $p=2$ \\
        \hline
        $N$ & $\frac{ L_2\Theta_1 \ln n}{\varepsilon}$ & $\frac{n L_2\Theta_2}{\varepsilon}$ \\
        \hline
        $m$ & $\max\left\{1,\frac{\sigma^2}{\varepsilon L_2}\right\}$ & $\max\left\{1,\frac{\sigma^2}{\varepsilon L_2}\right\}$\\
        \hline
        $\Delta_\zeta$ & $\min\left\{\frac{(\ln n)^2}{n}L_2^2\Theta_1,\, \frac{\varepsilon^2}{n\Theta_1},\, \frac{\varepsilon L_2}{n}\right\}$ & $\min\left\{nL_2^2\Theta_2,\, \frac{\varepsilon^2}{n\Theta_2},\, \frac{\varepsilon L_2}{n}\right\}$\\
        \hline
        $\Delta_\eta$ & $\min\left\{\frac{\ln n}{\sqrt{n}}L_2\sqrt{\Theta_1},\, \frac{\varepsilon}{\sqrt{n\Theta_1}},\, \sqrt{\frac{\varepsilon L_2}{n}}\right\}$ & $\min\left\{\sqrt{n}L_2\sqrt{\Theta_2},\, \frac{\varepsilon}{\sqrt{n\Theta_2}},\, \sqrt{\frac{\varepsilon L_2}{n}}\right\}$\\
        \hline
        $Nm$ & $\max\left\{\frac{L_2\Theta_1 \ln n}{\varepsilon}, \frac{\sigma^2\Theta_1\ln n}{\varepsilon^2}\right\}$ & $\max\left\{\frac{nL_2\Theta_2}{\varepsilon}, \frac{n\sigma^2\Theta_2 }{\varepsilon^2}\right\}$\\
        \hline
    \end{tabular}
    \caption{Algorithm~\ref{Alg:RDFDS} parameters for the cases $p=1$ and $p=2$.}
    \label{tab:RDD}
\end{table}
}

\subsubsection{Algorithms and main results for strongly convex problems.}
To obtain faster rates, we assume additionally that $f$ is $\mu_p$-strongly convex w.r.t. $p$-norm.
Our algorithms and proofs rely on the following assumption. Let $x_*$ be some fixed point and $x$ be a random point such that $\EE_x \big[ \| x-x_* \|_p^2 \big] \leqslant R_p^2$, then
\begin{equation}
    \EE_x  d\left( \frac{x-x_*}{R_p} \right)  \leqslant \frac{\Omega_p}{2},
		\label{eq:expdUpBound}
\end{equation}
where $\EE_x$ denotes the expectation with respect to random vector $x$ and $\Omega_p$ is defined as follows.
For $p=1$ and our choice of the prox-function \eqref{eq:dp1}, $\Omega_p = {\rm e}n^{(\kappa-1)(2-\kappa)/\kappa}\ln n = O(\ln n)$ with  $\kappa=1 + \frac{1}{\ln n}$, see \cite{nemirovsky1983problem,juditsky2014deterministic}.
For $p=2$ and our choice of the prox-function, $\Omega_p = 1$.
Our Accelerated Randomized Directional Derivative method for strongly convex problems (ARDDsc)  is listed as Algorithm~\ref{ACDS_sc}.
\begin{algorithm}
	\caption{Accelerated Randomized Directional Derivative method for strongly convex functions  (ARDDsc)}\label{ACDS_sc}
	\begin{algorithmic}[1]
		\REQUIRE $x_0$~---starting point s.t. $\| x_0 - x_* \|_p^2 \leq R_p^2$; $K \geqslant 1$~--- number of iterations; $\mu_p$ -- strong convexity parameter.
		\ENSURE point $u_K$.		
		\STATE Set $				N_0 = \left\lceil \sqrt{\frac{8 a L_2\Omega_p }{\mu_p}}\right\rceil$, 
			where $a=384n^2\rho_n$.
		\FOR{$k=0,\, \dots, \, K-1$}		
		\STATE $
							 m_k := \max \left\{1, \left\lceil \frac{32 \sigma^2 N_0 2^{k}}{nL_2\mu_p R_p^2 } \right\rceil \right\}, \quad R_k^2 := R_p^2 2^{-k} + \frac{4 \Delta}{\mu_p} \left(1-2^{-k} \right),
						$
		\STATE Set $d_k(x) = R_k^2d\left(\frac{x-u_k}{R_k}\right)$.
		\STATE Run ARDD with starting point $u_k$ and prox-function $d_k(x)$ for $N_0$ steps with batch size $m_k$. 
		\STATE Set $u_{k+1}=y_{N_0}$, $k=k+1$.		
		\ENDFOR
		\RETURN $u_K$
	\end{algorithmic}
\end{algorithm}

\begin{Th}
\label{Th:ACDS_sc_rate}
Let $f$ in problem \eqref{eq:PrSt} be $\mu_p$-strongly convex and ARDDsc method be applied to solve this problem. Then
    \begin{equation}\label{eq:ARDFDSSConv}
    \begin{array}{rl}
    \EE f(u_K) - f^* \leqslant 
    \frac{\mu_p  R_p^2}{2} \cdot 2^{-K} + 2 \Delta .
    \end{array}
    \end{equation}
    where 
    $\Delta = \frac{61N_0}{24L_2}\Delta_\zeta + \frac{122N_0}{3L_2}\Delta_\eta^2
        + \frac{12\sqrt{2nR_p^2\Omega_p}}{N_0^2}\left(\frac{\sqrt{\Delta_\zeta}}{2}+ 2\Delta_\eta\right)
        + \frac{N_0^2}{12n\rho_nL_2} \left(\frac{\sqrt{\Delta_\zeta}}{2} 
        + 2\Delta_\eta\right)^2$.
    Moreover, under an appropriate choice of $\Delta_\zeta$ and $\Delta_\eta$ s.t. $2 \Delta \leqslant \e/2$, the oracle complexity to achieve $\e$-accuracy of the solution is
    $$
\O\left(\max\left\{n^{\frac12+\frac{1}{q}}\sqrt{\frac{L_2\Omega_p }{\mu_p}}\log_2 \frac{\mu_p R_p^2 }{ \e},\frac{n^{\frac{2}{q}}\sigma^2 \Omega_p}{\mu_p \e}\right\}\right).
    $$
\end{Th}

Before we proceed to the non-accelerated method, we give the appropriate choice of the accuracy of the directional derivative evaluation $\Delta_\zeta$, $\Delta_\eta$ for ARDDsc to achieve an accuracy $\e$ of the solution. These values are chosen such that the r.h.s. of \eqref{eq:ARDFDSSConv} is smaller than $\e$. For simplicity we omit numerical constants and summarize the obtained values of the algorithm parameters in Table 3 below. The last row represents the total number of oracle calls, that is, the number of directional derivative evaluations.
{
\renewcommand{\arraystretch}{2}
\begin{table}[h]
\label{tab:ACDS_sc}
    \hspace{-3em}
    \begin{tabular}{|c|c|c|}
        \hline
         & $p=1$ & $p=2$ \\
        \hline
        $\Delta_\zeta$ & $\min\left\{\varepsilon\sqrt{\frac{L_2\mu_1}{n\ln n\Omega_1}},\, \varepsilon^2\frac{n(\ln n)^2 L_2^2\Omega_1}{R_1^2\mu_1^2},\, \varepsilon\cdot\frac{\mu_1}{n\Omega_1}\right\}$ & $\min\left\{\varepsilon\sqrt{\frac{L_2\mu_2}{n^2\Omega_2}},\, \varepsilon^2\frac{n^3 L_2^2\Omega_2}{R_2^2\mu_2^2},\, \varepsilon\cdot\frac{\mu_2}{n\Omega_2}\right\}$\\
        \hline
        $\Delta_\eta$ & $\min\left\{\sqrt{\varepsilon}\sqrt[4]{\frac{L_2\mu_1}{n\ln n\Omega_1}},\, \varepsilon\frac{\sqrt{n}\ln n L_2\sqrt{\Omega_1}}{R_1\mu_1},\, \sqrt{\varepsilon}\cdot\sqrt{\frac{\mu_1}{n\Omega_1}}\right\}$ & $\min\left\{\sqrt{\varepsilon}\sqrt[4]{\frac{L_2\mu_2}{n^2\Omega_2}},\, \varepsilon\frac{\sqrt{n^3} L_2\sqrt{\Omega_2}}{R_2\mu_2},\, \sqrt{\varepsilon}\cdot\sqrt{\frac{\mu_2}{n\Omega_2}}\right\}$\\
        \hline
        Calls & $\max\left\{\sqrt{\frac{n \ln n L_2\Omega_1 }{\mu_1}}\log_2 \frac{\mu_1 R_1^2 }{ \e},\frac{\sigma^2 \Omega_1 \ln n}{\mu_1 \e}\right\}$ & $\max\left\{n\sqrt{\frac{L_2\Omega_2 }{\mu_2}}\log_2 \frac{\mu_2 R_2^2 }{ \e},\frac{n\sigma^2 \Omega_2}{\mu_2 \e}\right\}$\\
        \hline
    \end{tabular}
    \caption{Algorithm~\ref{ACDS_sc} parameters for the cases $p=1$ and $p=2$.}
\end{table}
}

Our Randomized Directional Derivative method for strongly convex problems (RDDsc)  is listed as Algorithm~\ref{CDS_sc}.
\begin{algorithm}
	\caption{Randomized Directional Derivative method for strongly convex functions  (RDDsc)}\label{CDS_sc}
	\begin{algorithmic}[1]
		\REQUIRE $x_0$~---starting point s.t. $\| x_0 - x_* \|_p^2 \leq R_p^2$; $K \geqslant 1$~--- number of iterations; $\mu_p$ -- strong convexity parameter.
		\ENSURE point $u_K$.		
		\STATE Set $				N_0 = \left\lceil \frac{8 a L_2\Omega_p }{\mu_p}\right\rceil$, 
		where $a = 384n\rho_n$.
		\FOR{$k=0,\, \dots, \, K-1$}		
		\STATE $m_k := \max \left\{1, \left\lceil \frac{16 \sigma^2  2^{k}}{L_2\mu_p R_p^2 } \right\rceil \right\}, \quad R_k^2 := R_p^2 2^{-k} + \frac{4 \Delta}{\mu_p} \left(1-2^{-k} \right),
						$
		\STATE Set $d_k(x) = R_k^2d\left(\frac{x-u_k}{R_k}\right)$.
		\STATE Run RDD with starting point $u_k$ and prox-function $d_k(x)$ for $N_0$ steps with batch size $m_k$. 
		\STATE Set $u_{k+1}=y_{N_0}$, $k=k+1$.		
		\ENDFOR
		\RETURN $u_K$
	\end{algorithmic}
\end{algorithm}

\begin{Th}
\label{Th:CDS_sc_rate}
Let $f$ in problem \eqref{eq:PrSt} be $\mu_p$-strongly convex and RDDsc method be applied to solve this problem. Then
    \begin{equation}\label{eq:RDFDSSConv}
    \begin{array}{rl}
    \EE f(u_K) - f^* \leqslant 
    \frac{\mu_p  R_p^2}{2} \cdot 2^{-K} + 2 \Delta .
    \end{array}
    \end{equation}
    where $\Delta = \frac{n}{12L_2}\Delta_\zeta + \frac{4n}{3L_2}\Delta_\eta^2 + \frac{8\sqrt{2nR_p^2\Omega_p}}{N_0}\left(\frac{\sqrt{\Delta_\zeta}}{2} + 2\Delta_\eta\right) + \frac{N_0}{3L_2\rho_n}\left(\frac{\sqrt{\Delta_\zeta}}{2} + 2\Delta_\eta\right)^2$.
    Moreover, under an appropriate choice of $\Delta_\zeta$ and $\Delta_\eta$ s.t. $2 \Delta \leqslant \e/2$, the oracle complexity to achieve $\e$-accuracy of the solution is
    $$
\O\left(\max\left\{\frac{n^{\frac{2}{q}}L_2\Omega_p }{\mu_p}\log_2 \frac{\mu_p R_p^2 }{ \e},\frac{n^{\frac{2}{q}}\sigma^2 \Omega_p}{\mu_p \e}\right\}\right).
    $$
\end{Th}

Before we proceed, we give the appropriate choice of the accuracy of the directional derivative evaluation $\Delta_\zeta$, $\Delta_\eta$ for RDDsc to achieve an accuracy $\e$ of the solution. These values are chosen such that the r.h.s. of \eqref{eq:RDFDSSConv} is smaller than $\e$. For simplicity we omit numerical constants and summarize the obtained values of the algorithm parameters in Table 4 below. The last row represents the total number of oracle calls, that is, the number of directional derivative evaluations.
{
\renewcommand{\arraystretch}{2}
\begin{table}[h]
\label{tab:CDS_sc}
    \centering
    \begin{tabular}{|c|c|c|}
        \hline
         & $p=1$ & $p=2$ \\
        \hline
        $\Delta_\zeta$ & $\min\left\{\frac{\varepsilon L_2}{n},\, \varepsilon^2\frac{(\ln n)^2L_2^2}{nR_1^2\mu_1^2},\, \varepsilon\frac{\mu_1}{n\Omega_1}\right\}$ & $\min\left\{\frac{\varepsilon L_2}{n},\, \varepsilon^2\frac{nL_2^2}{R_2^2\mu_2^2},\, \varepsilon\frac{\mu_2}{n\Omega_2}\right\}$ \\
        \hline
        $\Delta_\eta$ & $\min\left\{\sqrt{\frac{\varepsilon L_2}{n}},\, \varepsilon\frac{\ln nL_2}{\sqrt{n}R_1\mu_1},\, \sqrt{\varepsilon\frac{\mu_1}{n\Omega_1}}\right\}$ & $\min\left\{\sqrt{\frac{\varepsilon L_2}{n}},\, \varepsilon\frac{\sqrt{n}L_2}{R_2\mu_2},\, \sqrt{\varepsilon\frac{\mu_2}{n\Omega_2}}\right\}$ \\
        \hline
        Calls & $\max\left\{\frac{L_2\Omega_1 \ln n }{\mu_1}\log_2 \frac{\mu_1 R_1^2 }{ \e},\frac{\sigma^2 \Omega_1}{\mu_1 \e}\right\}$ & $\max\left\{\frac{nL_2\Omega_2 }{\mu_2}\log_2 \frac{\mu_2 R_2^2 }{ \e},\frac{n\sigma^2 \Omega_2}{\mu_2 \e}\right\}$\\
        \hline
    \end{tabular}
    \caption{Algorithm~\ref{CDS_sc} parameters for the cases $p=1$ and $p=2$.}
\end{table}
}

\section{Primal-dual methods}
In this section, we focus on the developed primal-dual first-order methods for convex problems with linear constraints.

\subsection{Primal-dual methods for solving infinite-dimensional games}
The results of this subsection are published in \cite{dvurechensky2015primal-dual}.
Consider two moving objects with dynamics given by the following equations:
\begin{align}
&\dot{x}(t) = A_x(t) x(t) + B(t) u(t),  \dot{y}(t) = A_y(t) y(t) + C(t) v(t), \notag \\
&(x(0),y(0)) = (x_0,y_0).
\label{Dynamics}
\end{align}
Here $x(t) \in \R^n$, $y(t) \in \R^m$ are the phase
vectors of these objects, $u(t)$ is the control of the
first object (pursuer), and $v(t)$ is the control of the
second object (evader).
Matrices $A_x(t), A_y(t), B(t)$,
and $C(t)$ are continuous and have appropriate sizes. The
system is considered on the time interval $[0,\theta]$.
Controls are restricted in the following way $u(t) \in P
\subseteq \R^p$, $v(t) \in Q \subseteq \R^q \quad \forall t
\in [0,\theta]$.
We assume that $P,Q$ are closed, convex
sets.

The goal of the pursuer is to minimize the value of the functional:
\begin{equation}
F(u,v) + \Phi(x(\theta),y(\theta)):=\int_0^\theta{\tilde F(\tau,u(\tau),v(\tau)) d \tau} + \Phi(x(\theta),y(\theta)).
\label{price}
\end{equation}
The goal of the evader is the opposite. We need to find an
optimal guaranteed result for each object, which leads to
the problem of finding the saddle point of the above
functional. We assume the following:
\begin{itemize}
\item
$u(\cdot) \in L^2([0,\theta],\R^p)$, and $v(\cdot) \in
L^2([0,\theta],\R^q)$ (for the notation simplification we denote $L^2([0,\theta],\R^p)$ by $L^2_p$ and $L^2([0,\theta],\R^q)$ by $L^2_q$),
\item
the saddle point in this class of strategies exists,
\item
the function $F(u,v)$ is upper semi-continuous in $v$ and
lower semi-continuous in $u$,
\item
$\Phi(x,y)$ is continuous.
\end{itemize}

Denote by $V_x(t,\tau)$ the transition matrix of the first
system in  \eqref{Dynamics}. It is the unique solution of
the following matrix Cauchy problem
$$
\frac{d V_x(t,\tau)}{dt}  = A_x(t) V_x(t,\tau), \quad t \geq \tau , \quad V_x(\tau,\tau) = E.
$$
Here $E$ is the identity matrix. If the matrix $A_x(t)$ is constant, then 

\noindent $V_x(t,\tau) = e^{(t-\tau)A}$.

If we solve the first differential equation in
 \eqref{Dynamics}, then we can express $x(\theta)$ as a result
of the application of the linear operator $\B:
L^2_p \to \R^n$:
\begin{equation}\label{BuDef}
x(\theta) = V_x(\theta,0) x_0 + \int_0^\theta  V_x(\theta,\tau) B(\tau) u(\tau) d \tau :=\tilde x_0 + \B u.
\end{equation}
Below, we will use the conjugate operator $\B^{\ast}$ for the
operator $\B$. Let us find it explicitly. Let $\mu$ be a
$n$-dimensional vector. Then
\begin{align}
&\la \mu , \B u \ra = \la \mu, \int_0^\theta  V_x(\theta,\tau) B(\tau) u(\tau) d \tau \ra = \int_0^\theta \la \mu, V_x(\theta,\tau) B(\tau) u(\tau) \ra d \tau = \notag \\
& =\int_0^\theta \la B^{T}(\tau)  V_x^{T}(\theta,\tau) \mu ,u(\tau)  \ra d \tau = \la \B^{\ast} \mu, u \ra.  \notag
\end{align}
Note that the vector $\zeta (t) = V_x^{T}(\theta,t) \mu$ is the
solution of the following Cauchy problem:
$$
\dot{\zeta} (t) = - A_x^{T}(t) \zeta(t), \quad \zeta (\theta) = \mu , \quad t\in [0,\theta].
$$
So we can solve this ODE and find $\B^{\ast} \mu$ using
the obtained solution $\zeta(t)$ as $ \B^{\ast} \mu (t) =
B^{T}(t) \zeta(t)$.

In the same way, we introduce the transition matrix
$V_y(t,\tau)$ of the second system in  \eqref{Dynamics}, the
operator $\C: L^2_q \to \R^m$ defined by
the formula 

\noindent $\C v := \int_0^\theta  V_y(\theta,\tau) C(\tau)
v(\tau) d \tau$, and the vector $\tilde y_0 := V_y(\theta,0)
y_0$. The adjoint operator $\C^{\ast}$ also can be computed
using the solution of some ODE.

So below we study differential game problem in the following form:
\begin{equation}
\min_{u \in \U} \left[ \max_{v \in \V} \left\{ F(u,v) + \Phi(x,y) :y = \tilde y_0 + \C v\right\} :  x = \tilde x_0 + \B u \right],
\label{mainEquation}
\end{equation}
where 
$$
\U := \{ u(\cdot) \in L^2_p: u(t) \in P \quad \forall t \in [0,\theta]\}, \V := \{ v(\cdot) \in  L^2_q: v(t) \in Q \quad \forall t \in [0,\theta]\}
$$ are sets of admissible strategies of the players and $u \in \U$, $v \in \V$  mean $u(\cdot) \in \U$, $v(\cdot) \in \V$.
Our goal is to introduce a computational
method for finding an approximate solution of the problem  \eqref{mainEquation}.

First, we consider the problem
 \eqref{mainEquation} under two assumptions.

{\bf A1} The sets $P$ and $Q$ are bounded.

{\bf A2} In  \eqref{price} the functional $F(\cdot,v)$ is
convex for any fixed $v$, $F(u,\cdot)$ is concave for any
fixed $u$, $\Phi(\cdot,y)$ is convex for any fixed $y$,
and $\Phi(x,\cdot)$ is concave for any fixed $x$.

From {\bf A1}, since the norms of the operators $\B,\C$ are bounded, $x(\theta),y(\theta)$ are also bounded and we can equivalently reformulate the problem  \eqref{mainEquation} in the following way:
\begin{align}
&\min_{u \in \U, x \in X} \left[ \max_{v \in \V, y \in Y} \left\{ F(u,v) + \Phi(x,y) :y = \tilde y_0 + \C v\right\} :  x = \tilde x_0 + \B u \right]=\notag \\
&\max_{v \in \V, y \in Y}\left[ \min_{u \in \U, x \in X}  \left\{ F(u,v) + \Phi(x,y) :x = \tilde x_0 + \B u\right\} : y = \tilde y_0 + \C v  \right],
\label{initialProblemC}
\end{align}
where the sets $X$ and $Y$ are closed, convex and bounded.
Let us introduce the spaces of dual variables $\lambda \in \R^m$ and $\mu
\in \R^n$ corresponding to the linear constraints in the problem  \eqref{initialProblemC}, and some norms $\|\cdot\|_{\lambda}$ and $\|\cdot\|_{\mu}$ in these spaces. We define the norms in the dual space in the standard way
$$
\|s_{\lambda}\|_{\lambda,\ast} := \max \{ \la s_{\lambda}, \lambda\ra : \| \lambda \|_{\lambda} \leq 1 \}, \quad
\|s_{\mu}\|_{\mu,\ast} := \max \{ \la s_{\mu}, \mu\ra : \| \mu \|_{\mu} \leq 1 \}.
$$
In the simple case both the primal and the dual norm are Euclidean.

\begin{Lm} \label{Lm:adjointConv}
Let the Assumptions {\bf A1}, {\bf A2} hold. Also assume that the function $F(u,v)$ is upper semi-continuous in $v$ and
lower semi-continuous in $u$, the function $\Phi(x,y)$ is continuous, and that the sets $P$ and $Q$ are convex and closed.
Then the
problem  \eqref{initialProblemC} is equivalent to the
problem
\begin{equation}
\begin{array}{rl}
&\min_{\lambda} \max_{\mu} \{  \min_{u \in \U} \max_{v \in
\V}\left[ F(u,v) - \la \mu, \B u \ra + \la \lambda, \C v
\ra \right] \\
\\
& +  \min_{x \in X} \max_{y \in Y}\left[ \Phi(x,y) + \la
\mu, x \ra - \la \lambda, y \ra \right] - \la \mu, \tilde
x_0 \ra + \la \lambda,  \tilde y_0 \ra \},
\label{adjointProblemC}
\end{array}
\end{equation}
which we call the conjugate problem to
 \eqref{initialProblemC}.
\end{Lm}

We assume that the problems
\begin{align}
&\psi_1(\lambda,\mu) := \min_{u \in \U} \max_{v \in \V}\left[ F(u,v) - \la \mu, \B u \ra + \la \lambda, \C v \ra \right], \label{psi1DefC}\\
&\psi_2(\lambda,\mu) := \min_{x \in X} \max_{y \in Y}\left[ \Phi(x,y) + \la \mu, x \ra - \la \lambda, y \ra \right]\label{psi2DefC}
\end{align}
are rather simple so that they can be solved efficiently
or in a closed-form.
Note that the conjugate problem is
finite-dimensional
and the saddle point in the problems  \eqref{psi1DefC},
 \eqref{psi2DefC} exists for all $\lambda \in \R^m, \mu
\in \R^n$.

Note that the problem  \eqref{psi1DefC} has the following form
\begin{align}
&\min_{u \in\U} \max_{v \in \V}\left[ \int_0^\theta \left\{ \tilde F(\tau,u(\tau),v(\tau)) - \la \B^{\ast} \mu (\tau), u(\tau) \ra + \la \C^{\ast} \lambda (\tau), v(\tau) \ra \right\}  d \tau  \right]= \notag \\
& = \int_0^\theta \left\{\min_{u \in \U} \max_{v \in \V} \left[ \tilde F(\tau,u(\tau),v(\tau)) - \la \B^{\ast} \mu (\tau), u(\tau) \ra + \la \C^{\ast} \lambda (\tau), v(\tau) \ra  \right]  d \tau \right\},
\label{pointwise}
\end{align}
and it can be solved pointwise.


\subsubsection{Algorithm for convex-concave problem}
We assume that we are given some prox-function $d_{\lambda}(\lambda)$ with prox-center $\lambda_0$, which is strongly convex with convexity parameter $\sigma_{\lambda}$ in the given norm $\|\cdot\|_{\lambda}$.
For $\mu$ we introduce the similar assumptions.
Since $(\lambda^{\ast},\mu^{\ast})$ is the saddle point,
$(\lambda^{\ast},\mu^{\ast})$ is a weak solution to the following
variational inequality
$
\la g(\lambda, \mu), (\lambda - \lambda^{\ast},\mu -
\mu^{\ast}) \ra \geq 0, \quad \forall \lambda, \mu,
$
where $g(\lambda, \mu) := (\psi'_{\lambda}(\lambda, \mu), - \psi'_{\mu}(\lambda, \mu))$.
We apply the method of Simple Dual
Averages (SDA) from \cite{nesterov2009primal-dual} for finding an
approximate solution of the finite-dimensional problem
 \eqref{adjointProblemC}.
Let us choose some $\kappa \in ]0,1[$. We consider a space of $z
:=(\lambda,\mu)$ with the norm
\begin{equation}
\left\|z\right\|_z := \sqrt{\kappa \sigma_{\lambda} \left\| \lambda \right\|_{\lambda}^2+(1-\kappa) \sigma_{\mu} \left\| \mu \right\|_{\mu}^2},
\label{normZ}
\end{equation}
an oracle $g(z) := (g_{\lambda}(z),-g_{\mu}(z))$, a new prox-function\\ $d(z)~:=~\kappa d_{\lambda}(\lambda)~+~(1~-~\kappa)~ d_{\mu}(\mu)$, which is strongly convex with constant $\sigma_0=1$ with respect to the norm  \eqref{normZ}. We define $W := \R^m  \times  \R^n$.
The conjugate norm for  \eqref{normZ} is
$
\|g\|_{z,\ast} := \sqrt{\frac{1}{\kappa \sigma_{\lambda}} \|g_{\lambda}\|^2_{\lambda,\ast} + \frac{1}{(1-\kappa) \sigma_{\mu}} \|g_{\mu}\|^2_{\mu,\ast}}.
$
So we have a uniform upper bound for the answers of the
oracle $\left\|g(\lambda, \mu)\right\|_{z,\ast}^2~\leq~L^2~:=~\frac{L_{\lambda}^2}{\kappa \sigma_{\lambda}}+
\frac{L_{\mu}^2}{(1-\kappa) \sigma_{\mu} }$, where
$L_{\lambda} :=\sqrt{\theta} \left\| \C \right\|_{\lambda,L^2_q}\diam_2 Q  + \diam_{\lambda,\ast} Y + \| \tilde y_0 \|_{\lambda,\ast}$ and $L_{\mu} :=\sqrt{\theta} \left\| \B \right\|_{\mu,L^2_p} \diam_2 P  + \diam_{\mu,\ast} X + \| \tilde x_0 \|_{\mu,\ast}$.

The SDA method for solving  \eqref{adjointProblemC} is the following
\begin{enumerate}
    \item Initialization: Set $s_0 = 0$. Choose $z_0$, $\gamma > 0$.
    \item Iteration ($k \geq 0$):
    \subitem Compute $g_k = g(z_k)$. Set $s_{k+1} = s_k + g_k. \qquad \qquad \qquad \qquad \qquad {\rm (M1)}$
    \subitem $\beta_{k+1} = \gamma \hat{\beta}_{k+1}$. Set $z_{k+1} = \pi_{\beta_{k+1}}(-s_{k+1})$.
\end{enumerate}
Here the sequence $\hat{\beta}_{k+1}$ is defined by relations $\hat{\beta}_{0}=\hat{\beta}_{1}=1$, $\hat{\beta}_{i+1}=\hat{\beta}_{i}+\frac{1}{\hat{\beta}_{i}}$, for $i \geq 1$.
The mapping $\pi_{\beta}(s)$ is defined in the following
way
$
\pi_{\beta}(s) := \arg \min_{z \in W} \left\{ - \la s,z \ra + \beta d(z) \right\}.
$

Since the saddle point in the problem
 \eqref{initialProblemC} does exist, there exists a saddle
point $(\lambda^{\ast},\mu^{\ast})$ in the conjugate
problem  \eqref{adjointProblemC}. According to the Theorem 1 in
\cite{nesterov2009primal-dual}, the method ${\rm (M1)}$ generates a
bounded sequence $\{z_i\}_{i\geq 0}$. Hence, the sequences
$\{\lambda_i\}_{i\geq 0},\{\mu_i\}_{i\geq 0}$ are also
bounded. So we can choose $D_{\lambda}, D_{\mu}$ such that
$d_{\lambda}(\lambda_i)~\leq~D_{\lambda}$, $d_{\mu}(\mu_i)~\leq~D_{\mu}$ 
for all ${i\geq 0}$ and also, the pair
$(\lambda^{\ast},\mu^{\ast})$ is an interior solution:
$\Bb^{\lambda}_{r/\sqrt{\kappa \sigma_{\lambda}}}(\lambda^{\ast})~\subseteq~W_{\lambda}~:=~\left\{~
\lambda: ~d_{\lambda}(\lambda) ~\leq ~D_{\lambda}\right\}$,
and $\Bb^{\mu}_{r/\sqrt{(1-\kappa) \sigma_{\mu}}}(\mu^{\ast})~
\subseteq ~W_{\mu} ~:= ~\left\{ ~\mu:~
d_{\mu}(\mu) ~\leq ~D_{\mu} ~\right\}$ for some $r>0$. Then
we have $z^{\ast} ~:= ~(\lambda^{\ast},\mu^{\ast}) ~\in~
{\mathcal F}_D ~:=~\left\{~z\in W: ~d(z)~
\leq ~D ~\right\}$ with $D ~:= ~\kappa D_{\lambda}
+(1-\kappa)D_{\mu}$ and $\Bb^{z}_r(z^{\ast}) \subseteq {\mathcal
F}_D$.

Let us introduce a gap function
\begin{equation}
\delta_k(D) := \max_{z} \left\{ \sum_{i=0}^k{\la g_i,z_i-z
\ra : z \in {\mathcal F}_D} \right\}. \label{gapFunc}
\end{equation}
From the Theorem 2 in \cite{nesterov2009primal-dual}
we have
\begin{equation}
\frac{1}{k+1}\delta_k(D) \leq \frac{\hat{\beta}_{k+1}}{k+1}  \left( \gamma D + \frac{L^2}{2 \gamma} \right).
\label{deltakEstim}
\end{equation}

Denote
\begin{equation}
\left( \hat{u}_{k+1}, \hat{v}_{k+1}, \hat{x}_{k+1}, \hat{y}_{k+1}\right) :=
\frac{1}{k+1} \sum_{i=0}^k {\left( u_{i}, v_{i}, x_{i}, y_{i}\right)},
\label{hatPoints1}
\end{equation}
where $(u_i,v_i)$, $(x_i,y_i)$ are the saddle points at
the point $(\lambda_i,\mu_i)$ in  \eqref{psi1DefC} and
 \eqref{psi2DefC} respectively.
We define a function
\begin{equation}
\begin{array}{rl} & \phi(u,x,v,y) := \min_{\lambda} \max_{\mu} \{ 
F(u,v) + \Phi(x,y) + \la \mu, x - \tilde x_0 - \B u \ra +\\
&+ \la
\lambda, \C v + \tilde y_0 - y\ra :
 d_{\lambda}(\lambda) \leq D_{\lambda},d_{\mu}(\mu) \leq
D_{\mu} \}. \label{phiDef}
\end{array}
\end{equation}
%
Since $d_{\lambda}(\lambda^{\ast}) \leq D_{\lambda}$,
$d_{\mu}(\mu^{\ast}) \leq D_{\mu}$, and the conjugate
problem is equivalent to the initial one, we conclude that
the initial problem is equivalent to the problem
\begin{equation}
\min_{u\in \U,x \in X}\max_{v \in \V,y \in Y}\phi(u,x,v,y).
\label{auxilProblem}
\end{equation}

Let us introduce two auxiliary functions:
\begin{equation}
\xi(u,x) := \max_{v \in \V,y \in Y}\phi(u,x,v,y),
\label{xiDef}
\end{equation}
\begin{equation}
\eta(v,y) := \min_{u\in \U,x \in X}\phi(u,x,v,y).
\label{etaDef}
\end{equation}
Note that $\xi(u,x)$ is convex, $\eta(v,y)$ is concave,
and $\xi(u,x) ~\geq~
\phi(u^{\ast},x^{\ast},v^{\ast},y^{\ast}) ~\geq ~\eta(v,y)$ for all 
$u \in \U, v \in \V,x \in X,y \in Y$, where
$\phi(u^{\ast},x^{\ast},v^{\ast},y^{\ast})$ is the
solution to  \eqref{auxilProblem}.

\begin{Th} \label{convexTheorem}
Let the assumptions {\bf A1} and {\bf A2} be true. Then the points  \eqref{hatPoints1} generated by the method ${\rm (M1)}$ satisfy:
\begin{equation}
\xi(\hat{u}_{k+1},\hat{x}_{k+1}) - \eta(\hat{v}_{k+1},\hat{y}_{k+1}) \leq \frac{\hat{\beta}_{k+1}}{k+1} \left( \gamma D + \frac{L^2}{2 \gamma} \right),
\label{errorInFunc1}
\end{equation}
\begin{equation}
\begin{array}{l}
\left\| \tilde{x}_0 + \B \hat{u}_{k+1} -  \hat{x}_{k+1} \right\|_{\mu,\ast}
\leq \frac{\hat{\beta}_{k+1}\sqrt{\sigma_{\mu}}}{r(k+1)}
\left( \gamma D + \frac{L^2}{2  \gamma} \right),\\
\left\| \tilde{y}_0 + \C \hat{v}_{k+1} -  \hat{y}_{k+1} \right\|_{\lambda,\ast}
\leq \frac{\hat{\beta}_{k+1}\sqrt{\sigma_{\lambda}}}{r(k+1)}
\left( \gamma D + \frac{L^2}{2  \gamma} \right).
\label{errorInFeasibility1}
\end{array}
\end{equation}
\end{Th}

\subsubsection{Algorithm for strongly convex-concave problem}
In this subsection, we consider the problem \eqref{mainEquation}, under stronger assumptions and obtain faster convergence rates.

{\bf A3} The function $F(\cdot,v)$ is strongly convex for any
fixed $v$ with constant $\sigma_{F_u}$ which does not
depend on $v$, and function $F(u,\cdot)$ is strongly
concave for any fixed $u$ with constant $\sigma_{F_v}$
which does not depend on $u$. Assume that:
\begin{align}
& \left\| \nabla_{u} F(u,v_1)- \nabla_{u} F(u,v_2) \right\|_{L^2_p} \leq L_{uv} \left\|v_1-v_2 \right\|_{L^2_q}, \label{nablauFLips1}\\
& \left\| \nabla_{v} F(u_1,v)- \nabla_{v} F(u_2,v) \right\|_{L^2_q} \leq L_{vu} \left\|u_1-u_2 \right\|_{L^2_p}. \label{nablauFLips2}
\end{align}

{\bf A4} $\Phi(\cdot,y)$ is strongly convex for any fixed $y$ with respect to the norm $\|\cdot\|_{\mu,\ast}$ with constant $\sigma_{\Phi x}$ which doesn't depend on $y$ and $\Phi(x,\cdot)$ is strongly concave for any fixed $x$ with respect to the norm $\|\cdot\|_{\lambda,\ast}$ with constant $\sigma_{\Phi y}$ which doesn't depend on $x$.
Also we assume that:
\begin{align}
& \left\| \nabla_{x} \Phi(x,y_1)- \nabla_{x} \Phi(x,y_2) \right\|_\mu \leq L_{xy} \left\|y_1-y_2 \right\|_{\lambda,\ast}, \label{nablaxPhiLips1}\\
& \left\| \nabla_{y} \Phi(x_1,y)- \nabla_{y} \Phi(x_2,y) \right\|_\lambda \leq L_{yx} \left\|x_1-x_2 \right\|_{\mu,\ast}, \label{nablaxPhiLips2}
\end{align}
\begin{align}
& \left\| \nabla_{x} \Phi(x_1,y)- \nabla_{x} \Phi(x_2,y) \right\|_{\mu} \leq L_{xx} \left\|x_1-x_2 \right\|_{\mu,\ast}, \label{nablaxPhiLips3}\\
& \left\| \nabla_{y} \Phi(x,y_1)- \nabla_{y} \Phi(x,y_2) \right\|_{\lambda} \leq L_{yy} \left\|y_1-y_2 \right\|_{\lambda,\ast}. \label{nablaxPhiLips4}
\end{align}

Note that the assumptions {\bf A3}, {\bf A4} imply that the
level sets of the functions $F(u,v), \Phi(x,y)$ are closed,
convex and bounded. Similarly to the proof of the Lemma
\ref{Lm:adjointConv}, we get that the conjugate problem
for  \eqref{mainEquation} is
\begin{equation}
\begin{array}{rl}
\min_{\lambda} \max_{\mu} \{ & \min_{u \in \U} \max_{v \in
\V}\left[ F(u,v) - \la \mu, \B u \ra + \la \lambda, \C v
\ra \right] \\
& +  \min_{x} \max_{y}\left[ \Phi(x,y) + \la \mu, x \ra -
\la \lambda, y \ra \right] - \la \mu, \tilde x_0 \ra + \la
\lambda,  \tilde y_0 \ra \; \}. \label{adjointProblemSC}
\end{array}
\end{equation}
Here $\lambda \in \R^m$ and $\mu \in \R^n$.

We assume that the problems
\begin{align}
&\psi_1(\lambda,\mu) := \min_{u \in \U} \max_{v \in \V}\left[ F(u,v) - \la \mu, \B u \ra + \la \lambda, \C v \ra \right], \label{psi1DefSC}\\
&\psi_2(\lambda,\mu) := \min_{x} \max_{y}\left[ \Phi(x,y) + \la \mu, x \ra - \la \lambda, y \ra \right]\label{psi2DefSC}
\end{align}
are simple, which means that they can be solved
efficiently or in a closed form. Note that 
the saddle points in the problems
 \eqref{psi1DefSC},  \eqref{psi2DefSC} exists for all
$\lambda \in \R^m$ and $\mu \in \R^n$.

We assume that the norms $\|\cdot\|_{\lambda}$ and $\|\cdot\|_{\mu}$ are Euclidian.
Let us introduce the prox-function $d_{\lambda}(\lambda) :=
\frac{\sigma_{\lambda}}{2} \left\|\lambda \right\|_{\lambda}^2$. The function
$d_{\lambda}(\lambda)$ is strongly convex in this norm
with the convexity parameter $\sigma_{\lambda}$. For the
variable $\mu$ we introduce the prox-function $d_{\mu}(\mu) :=
\frac{\sigma_{\mu}}{2} \left\|\mu \right\|^2_{\mu}$, which is
strongly convex with the convexity parameter $\sigma_{\mu}$
with respect to the norm $\|\cdot\|_{\mu}$. These prox-functions are
differentiable everywhere.

For any $\lambda_1,\lambda_2 \in \R^m$ we can define the
Bregman distance:
$$
\omega_{\lambda}(\lambda_1,\lambda_2) := d_{\lambda}(\lambda_2)-d_{\lambda}(\lambda_1)-\la \nabla d_{\lambda}(\lambda_1),\lambda_2-\lambda_1\ra.
$$
Using the explicit expression for $d_{\lambda}(\lambda)$, we
get $\omega_{\lambda}(\lambda_1,\lambda_2)=
\frac{\sigma_{\lambda}}{2}\left\|
\lambda_1-\lambda_2\right\|^2$. Let us choose
$\bar{\lambda}=0$ as the center of the space $\R^m$. Then
we have
$\omega_{\lambda}(\bar{\lambda},\lambda)=d_{\lambda}(\lambda)$.
For $\mu$ we introduce the similar settings.

In the same way as it was done above, we conclude that finding the saddle
point $(\lambda^{\ast},\mu^{\ast})$ for the conjugate problem
 \eqref{adjointProblemSC} is equivalent to solving the
variational inequality
\begin{equation}
\la g(\lambda, \mu), (\lambda - \lambda^{\ast},\mu - \mu^{\ast}) \ra \geq 0, \quad \forall \lambda, \mu,
\label{varIneq2}
\end{equation}

\begin{equation}
\text{where} \;\;\; g(\lambda, \mu) := (\nabla_{\lambda} \psi(\lambda, \mu), - \nabla_{\mu} \psi(\lambda, \mu)).
\label{oracle2}
\end{equation}

Let us choose some $\kappa \in ]0,1[. $
Consider a space of $z :=(\lambda,\mu)$ with the norm
$$
\left\|z\right\|_z := \sqrt{\kappa \sigma_{\lambda} \left\| \lambda \right\|^2_{\lambda}+(1-\kappa) \sigma_{\mu} \left\| \mu \right\|^2_{\mu}},
$$
an oracle $g(z) := (\nabla_{\lambda} \psi(\lambda, \mu),- \nabla_{\mu} \psi(\lambda, \mu))$, a new prox-function
$$
d(z) := \kappa d_{\lambda}(\lambda) + (1 - \kappa) d_{\mu}(\mu)
$$
which is strongly convex with constant $\sigma_0 = 1$. We define $W := \R^m  \times \R^n$, the Bregman distance
$$
\omega(z_1,z_2):=\kappa \omega_{\lambda}(\lambda_1,\lambda_2) + (1 - \kappa) \omega_{\lambda}(\mu_2,\mu_2)
$$
which has an explicit form of $\omega(z_1,z_2)=d(z_1-z_2)$, and center $\bar{z}=(0,0)$. Then, $\omega(\bar{z},z)=d(z)$.
Note that the norm in the dual space is defined as
$$
\left\|g \right\|_{z,\ast} := \sqrt{\frac{1}{\kappa \sigma_{\lambda}} \left\| g_{\lambda} \right\|^2_{\lambda,\ast}+\frac{1}{(1-\kappa) \sigma_{\mu}} \left\| g_{\mu} \right\|^2_{\mu,\ast}}.
$$

In accordance to \cite{nesterov2007dual} for solving
 \eqref{varIneq2}, we can use the following method:
\begin{enumerate}
    \item Initialization: Fix $\beta = L$. Set $s_{-1} = 0$.
    \item Iteration ($k \geq 0$):
    \subitem Compute $x_k = T_{\beta} (\bar{z},s_{k-1}), \qquad \qquad \qquad \qquad \qquad \qquad {\rm (M2)}$
    \subitem Compute $z_k = T_{\beta} (x_k,-g(x_{k}))$,
    \subitem Set $s_k=s_{k-1}-g(z_k)$.
\end{enumerate}
Here
$
T_{\beta} (z,s) := \arg \max_{x \in W} \{\la s,x-z \ra -\beta \omega (z,x) \}.
$

Similarly to \cite{nesterov2009primal-dual}, we can prove that the
method ${\rm (M2)}$ generates a bounded sequence
$\{z_i\}_{i\geq 0}$. Hence the sequences
$\{\lambda_i\}_{i\geq 0},\{\mu_i\}_{i\geq 0}$ are also
bounded. Also, since the saddle point in the problem
 \eqref{mainEquation} exists, there exists a saddle point
$(\lambda^{\ast},\mu^{\ast})$ for the conjugate problem
 \eqref{adjointProblemSC}. These arguments allow us to
choose $D_{\lambda}, D_{\mu}$ such that
$d_{\lambda}(\lambda_i) \leq D_{\lambda}$, $d_{\mu}(\mu_i)
\leq D_{\mu}$ for all ${i\geq 0}$, which also ensure that
$(\lambda^{\ast},\mu^{\ast})$ is an interior solution:
$$
\begin{array}{l}
\Bb^{\lambda}_{r/\sqrt{\kappa \sigma_{\lambda}}}(\lambda^{\ast})
\subseteq W_{\lambda} := \left\{
\lambda: d_{\lambda}(\lambda) \leq D_{\lambda}\right\},\\
\\
\Bb^{\mu}_{r/\sqrt{(1-\kappa) \sigma_{\mu}}}(\mu^{\ast}) \subseteq
W_{\mu} := \left\{ \mu: d_{\mu}(\mu)
\leq D_{\mu} \right\}
\end{array}
$$
for some $r>0$. Then we have $z^{\ast} :=
(\lambda^{\ast},\mu^{\ast}) \in {\mathcal F}_D
:=\left\{z\in W: d(z) \leq D \right\}$
with $D := \kappa D_{\lambda} +(1-\kappa)D_{\mu}$ and
$\Bb^z_r(z^{\ast}) \subseteq {\mathcal F}_D$.

\begin{Th}
Let the Assumptions {\bf A3} and {\bf A4} be true, $\kappa
= \frac{\sigma_{\mu}}{\sigma_{\mu}+\sigma_{\lambda}}$, and 
\begin{equation}
\begin{array}{c}
L=\frac{\sigma_{\lambda}+\sigma_{\mu}}{\sigma_{\mu}\sigma_{\lambda}}
\sqrt{2\left(\frac{\left\| \C
\right\|^2_{\lambda,L^2_q}}{\sigma_{F_v}}+\frac{1}{\sigma_{\Phi_y}}+\frac{\left\|
\B \right\|_{\mu,L^2_p}\left\| \C \right\|_{\lambda,L^2_q} L_{vu}} {\sigma_{F_u}
\sigma_{F_v}}+ \frac{L_{yx}}{\sigma_{\Phi x}
\sigma_{\Phi_y}} \right) }  \\
\sqrt{\left( \frac{\left\| \B
\right\|_{\mu,L^2_p}\left\| \C \right\|_{\lambda,L^2_q} L_{uv}}{\sigma_{F_u}
\sigma_{F_v}} + \frac{L_{xy}}{\sigma_{\Phi_x} \sigma_{\Phi
y}} +  \frac{\left\| \B \right\|^2_{\mu,L^2_p}}{\sigma_{F_u}} +
\frac{1}{\sigma_{\Phi_x}} \right) }.
\label{gLipsitzConstant}
\end{array}
\end{equation}
Let the points
$z_i=(\lambda_i, \mu_i), i \geq 0$ be generated by the method
${\rm (M2)}$. Let the points in  \eqref{hatPoints1} be defined
by points $(u_i,v_i)$, $(x_i,y_i)$ which are the saddle
points at the points $(\lambda_i,\mu_i)$ in
 \eqref{psi1DefSC} and  \eqref{psi2DefSC} respectively. Then
for functions $\xi(u,x),\eta(v,y)$ defined in  \eqref{xiDef}
and  \eqref{etaDef} we have:
\begin{equation}
\xi(\hat{u}_{k+1},\hat{x}_{k+1}) - \eta(\hat{v}_{k+1},\hat{y}_{k+1}) \leq \frac{LD}{k+1}.
\label{errorInFunc2}
\end{equation}
Also the following is true:
\begin{equation*}
\left\| \B \hat{u}_{k+1} + \tilde x_0 -  \hat{x}_{k+1} \right\|_{\mu,\ast} \leq \frac{LD\sqrt{\sigma_{\mu}}}{r(k+1)}, \quad
\left\| \C \hat{v}_{k+1} + \tilde y_0 -  \hat{y}_{k+1} \right\|_{\lambda,\ast} \leq \frac{LD\sqrt{\sigma_{\lambda}}}{r(k+1)}.
\end{equation*}
\end{Th}

\subsection{Accelerated primal-dual gradient method for strongly convex problems with linear constraints}
The results of this subsection are published in \cite{chernov2016fast,dvurechensky2018computational}. See also a close work \cite{dvurechensky2016primal-dual}.

The main motivation for the algorithms in this subsection is approximating the optimal transport (OT) distance, which amounts to solving the \emph{OT problem} \cite{kantorovich1942translocation}:
\begin{align}
\label{eq:OT}
\centering
& \quad \quad \quad \quad \quad \min_{X \in \U(r,c)} \la C, X \ra, \notag \\
&\U(r,c) := \{X \in \R^{n \times n}_+ : \, X\one = r, \, X^T\one = c\},
\end{align}
where $X$ is \emph{transportation plan}, $C \in \R^{n \times n}_+$ is a given ground cost matrix, $r,c \in \R^n$ are given vectors from the probability simplex $\Delta^{n}$, $\one$ is the vector of all ones. 
The \emph{regularized OT problem} is
\begin{align}
\label{eq:ROT}
& \min_{X \in \U(r,c)} \la C, X \ra + \gamma \Reg(X), 
\end{align}
where $\gamma >0$ is the \emph{regularization parameter} and $\Reg(X)$ is a strongly convex \emph{regularizer}, e.g. negative entropy or squared Euclidean norm.
Our goal is to find $\widehat{X} \in \U(r,c)$ such that
\begin{equation}
\label{eq:OTSolDef}
\la C, \widehat{X} \ra \leq \min_{X \in \U(r,c)} \la C, X \ra + \e.
\end{equation}
In this case, $\la C, \widehat{X} \ra$ is an $\e$-approximation for the OT distance and $\widehat{X}$ is an approximation for the transportation plan.

Let us introduce some notation.For a general finite-dimensional real vector space $E$, we denote by $E^*$ its dual, given by linear pairing $\la g,x\ra$, $x\in E$, $g \in E^*$; by $\|\cdot\|_E$ the norm in $E$ and by $\|\cdot\|_{E,*}$ the norm in $E^*$, which is dual to $\|\cdot\|_E$. 
For a linear operator $A: E \to H$, we define its norm as 
$\|A\|_{E \to H} = \max_{x \in E,u \in H^*} \{\la u, A x \ra : \|x\|_{E} = 1, \|u\|_{H,*} = 1 \}$.
We say that a function $f: E \to \R$ is $\gamma$-strongly convex on a set $Q \subseteq E$ w.r.t. a norm in $E$ iff, for any $x,y \in Q$, $f(y) \geq f(x) + \la \nabla f(x) , y-x \ra + \frac{\gamma}{2}\|x-y\|^2_E$, where $\nabla f(x)$ is any subgradient of $f(x)$ at $x$.

For a matrix $A$ and a vector $a$, we denote $e^A$, $e^a$, $\ln A$, $\ln a$ their entrywise exponents and natural logarithms respectively. For a vector $a \in  \R^n$, we denote by $\|a\|_1$ the sum of absolute values of its elements, and by $\|a\|_2$ its Euclidean norm, and by $\|a\|_{\infty}$ the maximum absolute value of its elements. Given a matrix $A \in \R^{n \times n}$, we denote by ${\rm vec}(A)$ the vector in $\R^{n^2}$, which is obtained from $A$ by writing its columns one below another. For a matrix $A \in \R^{n \times n}$, we denote $\|A\|_1 = \|{\rm vec}(A)\|_1$ and $\|A\|_{\infty}= \|{\rm vec}(A)\|_{\infty}$. Further, we define the entropy of a matrix $X \in \R_+^{n \times n}$ by
\vspace{-0.8em}
\begin{equation}
\label{eq:Entr}
H(X) := - \sum_{i, j = 1}^n X^{i j} \ln X^{i j}.
\end{equation}
For two matrices $A,B$, we denote their Frobenius inner product by $\la A, B \ra$.
We denote by $\Delta^{n}: = \{a \in \R^n_+: a^T \one = 1\}$ the probability simplex in $\R^n$.


We start by considering a general minimization problem of strongly convex objective with linear constraints
\vspace{-1mm}
\begin{equation}
\label{eq:PrStGen}
\min_{x \in Q \subseteq E} \left\{ f(x) : Ax = b \right\},
\end{equation}
where $E$ is a finite-dimensional real vector space, $Q$ is a simple closed convex set, $A$ is a given linear operator from $E$ to some finite-dimensional real vector space $H$, $b \in H$ is given, $f(x)$ is a $\gamma$-strongly convex function on $Q$ with respect to some chosen norm $\|\cdot\|_E$ on $E$. 
\newline
The Lagrange dual problem for \eqref{eq:PrStGen}, written as a minimization problem, is
\vspace{-1mm}
\begin{equation}
\label{eq:DualPr}
\min_{\lambda \in H^*} \left\{ \vp(\lambda) := \la \lambda, b \ra + \max_{x\in Q} \left(- f(x) - \la A^T \lambda,x \ra \right) \right\}.
\end{equation}
Note that $\nabla \vp(\lambda) = b - Ax(\lambda)$ is Lipschitz-continuous \cite{nesterov2005smooth} 
\vspace{-1mm}
\[
\|\nabla \vp(\lambda_1) - \nabla \vp(\lambda_1) \|_H \leq L \|\lambda_1 - \lambda_2 \|_{H,*},
\]
where $x(\lambda) := \arg \min_{x \in Q} \left(- f(x) - \la A^T \lambda,x \ra \right)$ and $ L \leq \frac{\|A\|_{E \to H}^2}{\gamma}$. This estimate can be pessimistic and our algorithm does not use it and adapts automatically to the local value of the Lipschitz constant. 

We assume that the dual problem \eqref{eq:DualPr} has a solution and there exists some $R>0$ such that $\|\lambda^{*}\|_{2} \leq R < +\infty$,
where $\lambda^{*}$ is the solution to \eqref{eq:DualPr} with minimum value of $\|\lambda^{*}\|_{2}$. 
Note that the algorithm does not need any estimate of $R$ and the value $R$ is used only in the convergence analysis.

\begin{algorithm}[tb]
\caption{Adaptive Primal-Dual Accelerated Gradient Descent (APDAGD)}
\label{Alg:PDASTM}
\begin{algorithmic}[1]
   \REQUIRE Accuracy ${\e}_f,{\e}_{eq}> 0$, initial estimate $L_0$ s.t. $0<L_0<2L$.
   \STATE Set $i_0=k=0$, $M_{-1}=L_0$, $\beta_0=\alpha_0=0$, $\eta_0=\zeta_0=\lambda_0=0$.
   \REPEAT[Main iterate]
			\REPEAT[Line search]
				\STATE Set $M_k=2^{i_k-1}M_k$, find $\alpha_{k+1}$ s.t. $\beta_{k+1}:=\beta_k+\alpha_{k+1} = M_k\alpha_{k+1}^2$.
				Set $\tau_k = \alpha_{k+1}/\beta_{k+1}$.
				\STATE $\lambda_{k+1} = \tau_k\zeta_k + (1-\tau_k)\eta_k$. 
				\STATE $\zeta_{k+1} = \zeta_{k} - \alpha_{k+1} \nabla \vp(\lambda_{k+1})$.
				\STATE $\eta_{k+1} = \tau_k\zeta_{k+1} + (1-\tau_k)\eta_k$.
			\UNTIL{
				\begin{align}
					\vp(\eta_{k+1}) \leq & \vp(\lambda_{k+1}) + \la \nabla \vp(\lambda_{k+1}) ,\eta_{k+1} - \lambda_{k+1} \ra  + \frac{M_k}{2}\|\eta_{k+1} - \lambda_{k+1}\|_2^2. \notag				
				\end{align}
			}
			\STATE $\hat{x}_{k+1} =  \tau_kx(\lambda_{k+1})+(1-\tau_k)\hat{x}_{k}$.
			\STATE Set $i_{k+1}=0$, $k=k+1$.
  \UNTIL{$f(\hat{x}_{k+1})+\vp(\eta_{k+1}) \leq {\e}_f$, $\|A\hat{x}_{k+1}-b\|_{2} \leq {\e}_{eq}$.}
	\ENSURE $\hat{x}_{k+1}$, $\eta_{k+1}$.	
\end{algorithmic}
\end{algorithm}
\begin{Th}
\label{Th:PDASTMConv}
Assume that the objective in the primal problem \eqref{eq:PrStGen} is $\gamma$-strongly convex and that the dual solution $\lambda^*$ satisfies $\|\lambda^*\|_2 \leq R$. Then, for $k\geq 1$, the points $\hat{x}_{k}$, $\eta_k$ in Algorithm \ref{Alg:PDASTM} satisfy
\begin{align}
 f(\hat{x}_k) - f^* & \leq f(\hat{x}_k) + \vp(\eta_k) \leq \frac{16\|A\|_{E \to H}^2R^2}{\gamma k^2}, \label{eq:PDASTMObjBound} \\
\|A \hat{x}_k - b \|_2 &\leq \frac{16\|A\|_{E \to H}^2R}{\gamma k^2}, \label{eq:PDASTMConstrBound} \\
\|\hat{x}_k-x^*\|_E & \leq  \frac{8}{k} \frac{ \|A\|_{E \to H}R}{\gamma} ,\label{eq:PDASTMPlanBound}
\end{align}
where $x^*$ and $f^*$ are respectively an optimal solution and the optimal value in \eqref{eq:PrStGen}. Moreover, the stopping criterion in step 11 is correctly defined.
\end{Th}

Now we apply the general method to derive a complexity estimate for finding $\widehat{X} \in \U(r,c)$ satisfying \eqref{eq:OTSolDef}. 
We use entropic regularization of problem \eqref{eq:OT} and consider the regularized problem \eqref{eq:ROT} with the regularizer $\Reg(X) = - H(X)$, where $H(X)$ is given in \eqref{eq:Entr}. We define $E = \R^{n^2}$, $\|\cdot\|_E = \|\cdot\|_1$, and variable $x = {\rm vec}(X) \in \R^{n^2}$ to be the vector obtained from a matrix $X$ by writing each column of $X$ below the previous column. Also we set $f(x) = \la C,X\ra - \gamma H(X)$, $Q=\Delta^{n^2}$, $b^T = (r^T,c^T)$ and $A:\R^{n^2}\to \R^{2n}$ defined by the identity $(A\,{\rm vec}(X))^T = ((X \one)^T,(X^T \one)^T)$. With this setting, we solve problem \eqref{eq:PrStGen} by our APDAGD. Let $\widehat{X}_k$ be defined by identity ${\rm vec}(\widehat{X}_k) = \hat{x}_k$, where $\hat{x}_k$ is generated by APDAGD. We also define $\widehat{X} \in \U(r,c)$ to be the projection of $\widehat{X}_k$ onto $\U(r,c)$ constructed by Algorithm 2 in \cite{altschuler2017near-linear}. The pseudocode of our procedure for approximating the OT distance is listed as Algorithm \ref{Alg:OTbyGD}.

\begin{algorithm}[tb]
   \caption{Approximate OT by APDAGD}
   \label{Alg:OTbyGD}
\begin{algorithmic}[1]
	\REQUIRE{Accuracy $\e$.}
	\STATE Set $\gamma = \frac{\e}{3\ln n}.$
   \FOR{$k=1,2,...$}
	\STATE Make step of APDAGD and calculate $\widehat{X}_k$ and $\eta_k$.
	\STATE Find $\widehat{X}$ as the projection of $\widehat{X}_k$ on $\U(r,c)$ by Algorithm 2 in \cite{altschuler2017near-linear}.
   \IF{$\la C,\widehat{X}-\widehat{X}_k\ra \leq \frac{\e}{6}$ and $f(\hat{x}_k)+\vp(\eta_k) \leq \frac{\e}{6}$}
   \STATE Return $\widehat{X}$.
	\ELSE
	\STATE $k = k+1$ and continue.
   \ENDIF
   \ENDFOR
\end{algorithmic}
\end{algorithm}

\begin{Th}
\label{Th:OTbyGDCompl}
Algorithm \ref{Alg:OTbyGD} outputs $\widehat{X} \in \U(r,c)$ satisfying \eqref{eq:OTSolDef} in 
\begin{equation}
\label{eq:OTbyGDCompl}
O\left(\min\left\{\frac{n^{9/4}\sqrt{R\|C\|_{\infty}\ln n}}{\e}, \frac{n^2R\|C\|_{\infty}\ln n}{\e^2}\right\}\right)
\end{equation}
arithmetic operations.
\end{Th}

\subsection{Distributed primal-dual accelerated stochastic gradient method}
The results of this subsection are published in \cite{dvurechensky2018decentralize}.

In this subsection we are motivated by regularized semi-discrete formulation of the optimal transport problem. We start with some notation. We define $\M(\X)$ -- the set of positive Radon probability measures on a metric space $\X$, and $S_1(n)  = \{ a \in \mathbb{R}_+^n  \mid \sum_{l=1}^n a_l =1 \}$ the probability simplex. We use $\C(\X)$ as the space of continuous functions on $\X$. We denote by $\delta(x)$ the Dirac measure at point $x$. We refer to $\lambda_{\max}(W)$ as the maximum eigenvalue of matrix W. We also use bold symbols for stacked vectors $\p = [p_1^T,\cdots,p_m^T]^T \in \mathbb{R}^{mn}$, where $p_1,...,p_m\in \R^n$. In this case $[\p]_i=p_i$ -- the $i$-th block of $\p$. For a vector $\lambda \in \R^n$, we denote by $[\lambda]_l$ its $l$-th component. We refer to the Euclidean norm of a vector $\|p\|_2:=\sqrt{\sum_{l=1}^n([p]_l)^2}$ as $2$-norm. 

Following the line of work started by \cite{cuturi2013sinkhorn}, we consider entropic regularization for the optimal transport problem and the corresponding regularized Wasserstein distance and barycenter. Assume that we are given a positive Radon probability measure  $\mu$  with density $q(y)$  on a metric space $\mathcal{Y}$, and a discrete probability measure $\nu=\sum_{i=1}^np_i\delta(z_i)$ with weights $p$ and finite support given by points $z_1, \dots, z_n \in \mathcal{Z}$  from a metric space $\mathcal{Z}$. 
The regularized Wasserstein distance in semi-discrete setting between continuous measure $\mu$ and discrete measure $\nu$ is defined as 
\begin{align*}
\W_{\gamma}(\mu,\nu)=\min_{\pi \in \Pi(\mu,\nu)}\left\{\sum_{i=1}^n\int_{\Y}  c_i(y)\pi_i(y)dy   + \gamma KL(\pi|\xi)\right\},
\end{align*}
where $c_i(y) = c(z_i,y)$ is a cost function for transportation of a unit of mass from point $z_i$ to point $y$, $\xi$ is the uniform distribution on $\Y \times \mathcal{Z}$, 
{\sloppy $KL(\pi|\xi)=\sum_{i=1}^n\int_{\Y}  \pi_i(y)\log\left(\frac{\pi_i(y)}{\xi}\right)dy$}, and the set of admissible coupling measures $\pi$ is defined as follows
{\small \begin{align*}
\Pi(\mu,\nu) = \left\{\pi \in\M(\Y) \times  S_1(n): \sum_{i=1}^n \pi_i(y) = q(y), y \in \mathcal{Y}, \int_{\mathcal{Y}} \pi_i(y)dy = p_i, \forall~i=1,\dots,n  \right\}.   
\end{align*}}
For a set of positive Radon probability measures $(\mu_1, \dots, \mu_m)$ the regularized Wasserstein barycenter in the semi-discrete setting is defined as the solution $p$ to the following convex optimization problem
\begin{align}\label{w_barycenter}
\min_{p \in S_1(n)} \sum\limits_{i=1}^{m}  \W_{\gamma,\mu_i}(p) = \min_{\substack{p_1=\cdots=p_m \\ p_1,\dots,p_m \in S_1(n)}} 
\sum\limits_{i=1}^{m} \W_{\gamma,\mu_i}(p_i),
\end{align}
where we fixed the support $z_1, \dots, z_n \in \mathcal{Z}$ of the barycenter $\nu$ and characterize it by the vector $p \in S_n(1)$, i.e., $\nu=\sum_{i=1}^np_i\delta(z_i)$ and $\W_{\gamma,\mu}(p) := \W_{\gamma}(\mu,\nu)$.

We now describe the distributed optimization setting for solving the second problem in \eqref{w_barycenter}.
We assume that each measure $\mu_i$ is held by an agent $i$ on a network and this agent can sample from this measure. We model such a network as a fixed \textit{connected undirected graph} \mbox{$\mathcal{G} = (V,E)$}, where 
$V$ is the set of $m$ nodes and $E$ is the set of edges. We assume that the graph $\mathcal{G}$ does not have self-loops. The network structure imposes information constraints, specifically, each node $i$ has access to $\mu_i$ only and a node can exchange information only with its immediate neighbors, i.e., a node $i$ can communicate with node $j$ if and only if $(i,j)\in E$. 

We represent the communication constraints imposed by the network by introducing a single equality constraint instead of the constraints $p_1=\cdots=p_m$ in \eqref{w_barycenter}.
To do so, we define the Laplacian matrix 
$\bar W{\in \mathbb{R}^{m\times m}}$ of the graph $\mathcal{G}$ such that a) $[\bar W]_{ij} = -1$ if $(i,j) \in E$, b) $[\bar W]_{ij} = \text{deg}(i)$ if $i=j$, c) $[\bar W]_{ij} = 0$ otherwise.
Here $\text{deg}(i)$ is the degree of the node $i$, i.e., the number of neighbors of the node. 
Finally, define the communication matrix (also referred to as an interaction matrix) 
by \mbox{$W := \bar W \otimes I_n$}. 


In this setting, $\sqrt{W}{\mathtt{p}} = 0$ if and only if {$p_1 = \cdots = p_m$}, where we defined stacked column vector $\mathtt{p} = [p_1^T,\cdots,p_m^T]^T \in \mathbb{R}^{mn}$. 
Using this fact, we equivalently rewrite problem~\eqref{w_barycenter} as the maximization problem with linear equality constraint
	\begin{align}\label{consensus_problem2}
	\max_{\substack{p_1,\dots, p_m \in S_1(n) \\ \sqrt{W} \mathtt{p}=0 }} ~ - \sum\limits_{i=1}^{m} \W_{\gamma, \mu_i}(p_i) .
	\end{align}
    
    Given that problem \eqref{consensus_problem2} is an optimization problem with linear constraints, we introduce a stacked vector of dual variables $\Blm = [\lambda_1^T,\cdots,\lambda_m^T]^T \in \R^{mn}$ for the constraints $\sqrt{W}\p=0$ in \eqref{consensus_problem2}.
Then, the Lagrangian dual problem for \eqref{consensus_problem2} is
\begin{align}\label{eqWB:DualPr}
\min_{\Blm \in \R^{mn}} ~\max_{p_1,\dots, p_m \in S_1(n) } ~ \left\lbrace \sum\limits_{i=1}^{m} \langle \lambda_i, [\sqrt{W}\mathtt{p}]_i\rangle-\W_{\gamma, \mu_i}(p_i)\right\rbrace = \min_{\Blm \in \R^{mn}} \sum_{i=1}^{m} \W^*_{\gamma, \mu_i}([\sqrt{W}\Blm]_i),
\end{align}
where $[\sqrt{W}\p]_i$ and $[\sqrt{W}\Blm]_i$ denote the $i$-th $n$-dimensional block of vectors $\sqrt{W}\p$ and $\sqrt{W}\Blm$ respectively, and $\W^*_{\gamma,\mu_i}(\cdot)$ is the Fenchel-Legendre transform of $\W_{\gamma,\mu_i}(p_i)$.

Next, we consider a general smooth stochastic convex optimization problem which is dual to some optimization problem with linear equality constraints. 
For any finite-dimensional real vector space $E$, we denote by $E^*$ its dual. 
Let $\|\cdot\|_E$ denote some norm on $E$ and $\|\cdot\|_{E,*}$ denote the norm on $E^*$ which is dual to $\|\cdot\|_E$
$\|\lambda\|_{E,*} = \max_{\|x\|_E \leq 1} \la \lambda, x \ra$.
For a linear operator $A:E_1 \to E_2$, we define the adjoint operator $A^T: E_2^* \to E_1^*$ in the following way $\la u, A x \ra = \la A^T u, x \ra, \quad \forall u \in E_2^*, \quad x \in E_1$.
We say that a function $f: E \to \R$ has a $L$-Lipschitz-continuous gradient w.r.t. norm $\|\cdot\|_{E,*}$ if it is differentiable and its gradient satisfies Lipschitz condition $
\|\nabla f(x) - \nabla f(y) \|_{E,*} \leq L \|x-y\|_E, \quad \forall x,y \in E$.

Our next goal is to provide an algorithm for a primal-dual pair of problems
\begin{equation}
(P) \quad \quad \min_{x\in Q \subseteq E} \left\{ f(x) : Ax =b\right\}, \quad (D) \quad \min_{\lambda \in \Lambda} \left\{   \la \lambda, b \ra +  \max_{x\in Q} \left( -f(x) - \la A^T \lambda ,x \ra \right) \right\}.
\notag
\end{equation}
where $Q$ is a simple closed convex set, $A: E \to H$ is given linear operator, $b \in H$ is given, $\Lambda = H^*$.
We define
\begin{equation}
\vp(\lambda) :=\la \lambda, b \ra +  \max_{x\in Q} \left( -f(x) - \la A^T \lambda  ,x \ra \right) = \la \lambda, b \ra + f^*(-A^T\lambda)
\label{eq:vp_def}
\end{equation}
and assume it to be smooth with $L$-Lipschitz-continuous gradient. Here $f^*$ is the Fenchel-Legendre dual for $f$. 
We also assume that $f^*(-A^T\lambda)= \E_\xi F^*(-A^T\lambda,\xi)$, where $\xi$ is random vector. Also, we define $F(x,\xi)$ to be the Fenchel-Legendre conjugate function to $F^*$, i.e. it satisfies $F^*(-A^T\lambda,\xi) = \max_{x\in Q}\{\la -A^T \lambda,x\ra - F(x,\xi) \}$ and $x(\lambda,\xi)$ to be the solution of this maximization problem.
Under these assumptions, the dual problem $(D)$ can be accessed by a stochastic oracle $(\Phi(x,\xi),\nabla \Phi(\lambda,\xi)) = (F^*(-A^T\lambda,\xi),\nabla F^*(-A^T\lambda,\xi))$ satisfying $\E_\xi \Phi(\lambda,\xi) = \vp(\lambda)$, $\E_\xi \nabla \Phi(\lambda,\xi) = \nabla \vp(\lambda)$, which we use in our algorithm.
Finally, we assume that dual problem $(D)$ has a solution $\lambda^*$ and there exists some $R >0$ such that $\|\lambda^*\|_{2} \leq R < +\infty$.


We additionally assume that the variance of the stochastic approximation $\nabla \Phi(\lambda,\xi)$ for the gradient of $\vp$ can be controlled and made as small as we desire. This can be done, for example by mini-batching the stochastic approximation. Also, since $\nabla  \Phi(\lambda,\xi) =b - A \nabla F^*(-A^T\lambda,\xi) = b - A x(\lambda,\xi)$, on each iteration, to find $\nabla  \Phi(\lambda,\xi)$ we find the vector $x(\lambda,\xi)$ and use it for the primal iterates.

\begin{algorithm}[H]
\caption{Accelerated Primal-Dual Stochastic Gradient Method (APDSGM)}
\label{Alg:APDSGD}
\begin{algorithmic}[1]
   \REQUIRE Number of iterations $N$.
   \STATE $C_0=\alpha_0=0$, $\eta_0=\zeta_0=\lambda_0=0$.
                \FOR{$k=0,\dots, N-1$}
                \STATE Find $\alpha_{k+1}$ as the largest root of the equation $C_{k+1}:=C_k+\alpha_{k+1} = 2L\alpha_{k+1}^2$. 
                $\tau_{k+1}=\alpha_{k+1}/C_{k+1}$.
				\STATE
				$\lambda_{k+1} = \tau_{k+1}\zeta_k + (1-\tau_{k+1})\eta_k$
				\STATE 
				$\zeta_{k+1}= \zeta_{k} - \alpha_{k+1} \nabla  \Phi(\lambda_{k+1},\xi_{k+1})$.
				\STATE $
				\eta_{k+1} =\tau_{k+1}\zeta_{k+1} + (1-\tau_{k+1})\eta_k$.
			\STATE 
					$
                    \hat{x}_{k+1} = \tau_{k+1}x(\lambda_{k+1},\xi_{k+1})+(1-\tau_{k+1})\hat{x}_{k}$.
                \ENDFOR
	\ENSURE The points $\hat{x}_{k+1}$, $\eta_{k+1}$.	
\end{algorithmic}
\end{algorithm}

\begin{Th}
\label{Th:stoch_err}
Let $\vp$ have $L$-Lipschitz-continuous gradient w.r.t. 2-norm and $\|\lambda^*\|_2 \leq R$, where $\lambda^*$ is a solution of dual problem $(D)$. Given desired accuracy $\e$, assume that, at each iteration of Algorithm \ref{Alg:APDSGD}, the stochastic gradient $\nabla \Phi(\lambda_k,\xi_k)$ is chosen in such a way that $\E_\xi \|\nabla \Phi(\lambda_k,\xi_k) - \nabla \vp(\lambda_k) \|_2^2 \leq \frac{\e L\alpha_k}{C_k}$.
Then, for any $\varepsilon > 0$ and $N \geq 0$, and expectation $\E$ w.r.t. all the randomness $\xi_1,\dots, \xi_N$, the outputs $\eta_N$ and $\hat{x}_N$ generated by the  Algorithm \ref{Alg:APDSGD} satisfy  
\begin{align}\label{eq:APDSGDRate}
f(\mathbb{E}\hat{x}_N)-f^* \leq \frac{32LR^2}{N^2} + \frac{\e}{2}~ ~ ~ ~ \text{and} ~ ~ ~ ~ \|A\mathbb{E}\hat{x}_N-b\|_2 \leq \frac{32LR}{N^2} + \frac{\e}{2R},
\end{align}
\end{Th}
Next, we apply the general algorithm to solve the primal-dual pair of problems \eqref{consensus_problem2}-\eqref{eqWB:DualPr} and approximate the regularized Wasserstein barycenter which is a solution to \eqref{consensus_problem2}. 
\begin{Lm}\label{Lm:dual_obj_properties2}
The gradient of the objective function $\W_{\gamma}^*(\Blm)$ in the dual problem \eqref{eqWB:DualPr} is $\lambda_{\max}(W) /\gamma$-Lipschitz-continuous w.r.t. 2-norm. If its stochastic approximation is defined as
\begin{align} \label{eq:StochGradDef}
[\widetilde{\nabla} \W_{\gamma}^*(\Blm)]_i &=  \sum_{j=1}^{m}\sqrt{W}_{ij} \widetilde{\nabla} \W_{\gamma,\mu_j}^*(\blm_j), \; i=1,...,m, \; \text{with}\notag \\
\widetilde{\nabla} \W_{\gamma,\mu_j}^*(\blm_j) &= \frac{1}{M}\sum_{r=1}^{M} p_j(\blm_j), \; \text{and}\;
[p_j(\blm_j)]_l =
\frac{\exp(([\bar{\lambda}_j]_l-c_l(Y_r^j))/\gamma)  }{\sum_{\ell=1}^n\exp(([\bar{\lambda}_j]_{\ell}-c_\ell(Y_r^j))/\gamma)}.
\end{align}  
where $M$ is the batch size, $\blm_j := [\sqrt{W}\Blm]_j$, $j=1,...,m$, $Y_1^j,...,Y_r^j$ is a sample from the measure $\mu_j$, $j=1,...,m$. Then $\E_{Y_r^j\sim\mu_j, j=1,...,m, r=1,...,M} \widetilde{\nabla} \W_{\gamma}^*(\Blm) = \nabla \W_{\gamma}^*(\Blm)$ and
\begin{equation}\label{eq:variance}
\E_{Y_r^j\sim\mu_j, j=1,...,m, r=1,...,M} \|\widetilde{\nabla} \W_{\gamma}^*(\Blm) - \nabla \W_{\gamma}^*(\Blm)\|_2^2 \leq \frac{\lambda_{\max}(W)}{M}, \; \Blm \in \R^{mn}.
\end{equation}
\end{Lm}
Based on this lemma, we see that if, on each iteration of Algorithm \ref{Alg:APDSGD}, the mini-batch size $M_k$ satisfies $M_k \geq \frac{\lambda_{max}(W) C_{k}}{L\alpha_k\e}$, the assumptions of Theorem \ref{Th:stoch_err} hold.

For the particular problem \eqref{eqWB:DualPr} the step 5 of Algorithm \ref{Alg:APDSGD} can be written block-wise $[\Bzeta_{k+1}]_i= [\Bzeta_{k}]_i - \alpha_{k+1} \sum_{j=1}^{m}\sqrt{W}_{ij} \widetilde{\nabla} \W_{\gamma,\mu_j}^*([\sqrt{W}\Blm_{k+1}]_j)$, $i=1,...,m$. Unfortunately, this update can not be made in the decentralized setting since the sparsity pattern of $\sqrt{W}_{ij}$ can be different from $W_{ij}$ and this will require some agents to get information not only from their neighbors. To overcome this obstacle, we change the variables and denote $\bar{\Blm} = \sqrt{W}\Blm$, $\bar{\Beta} = \sqrt{W}\Beta$, $\bar{\Bzeta} = \sqrt{W}\Bzeta$. Then the step 5 of Algorithm \ref{Alg:APDSGD} becomes $[{\BBzeta}_{k+1}]_i= [{\BBzeta}_{k}]_i - \alpha_{k+1} \sum_{j=1}^{m}W_{ij} \widetilde{\nabla} \W_{\gamma,\mu_j}^*([\BBlm_{k+1}]_j)$, $i=1,...,m$.

\begin{algorithm}[H]
	\caption{Distributed computation of Wasserstein barycenter}
	\label{alg:main}
    {
	\begin{algorithmic}[1]
    \REQUIRE Each agent $i\in V$ is assigned its measure $\mu_i$.
		\STATE All agents set $[\BBeta_{0}]_i = [\BBzeta_{0}]_i = [\BBlm_{0}]_i = \boldsymbol{0} \in \mathbb{R}^n$, \\$C_0=\alpha_0=0$ and $N$ 
		\STATE{For each agent $i \in V$:}
		\FOR{ $k=0,\dots,N-1$ }
		\STATE Find $\alpha_{k+1}$ as the largest root of the equation \\$C_{k+1}:=C_k+\alpha_{k+1} = 2L\alpha_{k+1}^2$.\\ $\tau_{k+1}=\alpha_{k+1}/C_{k+1}$.        
		\STATE Set {$M_{k+1} = \max\left\{1,~{{\lambda_{\max}(W) C_{k+1}}/({L\alpha_{k+1}\e})}\right\}$}
        \STATE $[\BBlm_{k+1}]_i = \tau_{k+1}[\BBzeta_k]_i+(1-\tau_{k+1}) [\BBeta_k]_i$
		\STATE{Generate $M_{k+1}$ samples $\{Y_r^i\}_{r=1}^{M_{k+1}}$  from the measure $\mu_i$ and set $\widetilde{\nabla} \W_{\gamma,\mu_i}^*([\BBlm_{k+1}]_i)$ as in \eqref{eq:StochGradDef}.}
        \STATE{Share 
        $\widetilde{\nabla} \W_{\gamma,\mu_i}^*([\BBlm_{k+1}]_i)$ with $\{j \mid (i,j) \in E \}$}
		\STATE  $ [\BBzeta_{k+1}]_i = [\BBzeta_{k}]_i-\alpha_{k+1} \sum_{j=1}^{m} W_{ij} \widetilde{\nabla} \W_{\gamma,\mu_j}^*([\BBlm_{k+1}]_j)$
        \STATE  $[\BBeta_{k+1}]_i = \tau_{k+1}[\BBzeta_{k+1}]_i+(1-\tau_{k+1})[\BBeta_{k+1}]_i$
   \STATE $[\hat{\p}_{k+1}]_i  = \tau_{k+1}p_i([\BBlm_{k+1}]_i)+(1-\tau_{k+1})[\hat{\p}_{k+1}]_i$, where $p_i(\cdot)$ is defined in \eqref{eq:StochGradDef}.
   \ENDFOR
	\ENSURE $\hat{\p}_{N}$.
        
	\end{algorithmic}}
\end{algorithm}

\begin{Th}
\label{Th:WBCompl}
Under the above assumptions, Algorithm \ref{alg:main} after $N=\sqrt{16\lambda_{max}(W)R^2/(\e\gamma)}$ iterations returns an approximation $\hat{\p}_{N}$ for the barycenter, which satisfies
\begin{equation}\label{eq:WBerr}
\sum\limits_{i=1}^{m} \W_{\gamma,\mu_i}(\E[\hat{\p}_{N}]_i)-\sum\limits_{i=1}^{m} \W_{\gamma,\mu_i}([\p^*]_i) \leq \e, \quad \|\sqrt{W}\E\hat{\p}_{N}\|_2 \leq \e/R.
\end{equation}
Moreover, the total complexity is $O\left(n\max\eig R^2/\e^2, \sqrt{\eig R^2/(\e \gamma)}  \right)$ arithmetic operations.
\end{Th}

\subsection{Primal-dual accelerated gradient method with small-dimensional relaxation oracle}
The results of this subsection are published in \cite{guminov2019accelerated,nesterov2020primal-dual}. See also a follow-up work \cite{guminov2019acceleratedAM}.

In this subsection, we consider a minimization problem with linear equality constraints. 

Specifically, we consider the following minimization problem 
\begin{equation}
(P_1) \quad \quad \min_{x\in Q \subseteq E} \left\{ f(x) : \bm{A}x =b \right\},
\notag
\end{equation}
where $E$ is a finite-dimensional real vector space, $Q$ is a simple closed convex set, $\bm{A}$ is given linear operator from $E$ to some finite-dimensional real vector space $H$, $b \in H$ is given.
The Lagrange dual problem to Problem $(P_1)$ is
\begin{equation}
(D_1) \quad \quad \max_{\lambda \in \Lambda} \left\{ - \la \lambda, b \ra  + \min_{x\in Q} \left( f(x) + \la \bm{A}^T \lambda  ,x \ra \right) \right\}.
\notag
\end{equation}
Here we denote $\Lambda=H^*$.
It is convenient to rewrite Problem $(D_1)$ in the equivalent form of a minimization problem
\begin{align}
& (P_2) \quad \min_{\lambda \in \Lambda} \left\{   \la \lambda, b \ra  + \max_{x\in Q} \left( -f(x) - \la \bm{A}^T \lambda  ,x \ra \right) \right\}. \notag
\end{align}
We denote
\begin{equation}
\vp(\lambda) =  \la \lambda, b \ra  + \max_{x\in Q} \left( -f(x) - \la \bm{A}^T \lambda  ,x \ra \right).
\label{eqLS:vp_def}
\end{equation}
Since $f$ is convex, $\vp(\lambda)$ is a convex function and, by Danskin's theorem, its subgradient is equal to (see e.g. \cite{nesterov2005smooth})
\begin{equation}
\nabla \vp(\lambda) = b - \bm{A} x (\lambda)
\label{eq:nvp}
\end{equation}
where $x (\lambda)$ is some solution of the convex problem
\begin{equation}
\max_{x\in Q} \left( -f(x) - \la \bm{A}^T \lambda  ,x \ra \right).
\label{eq:inner}
\end{equation}

In what follows, we make the following assumptions about the dual problem $(D_1)$
\begin{itemize}
    \item Subgradient of the objective function $\vp(\lambda)$ satisfies H\"older condition with constant $M_{\nu}$, i.e., for all $\lambda, \mu \in \Lambda$ and some $\nu\in\left[0,1\right]$
\begin{equation}\label{eq:hold_cond}
\|\nabla \vp(\lambda) - \nabla \vp(\mu)\|_* \leqslant M_\nu \|\lambda-\mu\|^\nu.
\end{equation}
    \item The dual problem $(D_1)$ has a solution $\lambda^*$ and there exist some $R>0$ such that
	\begin{equation}
	\|\lambda^{*}\|_{2} \leqslant R < +\infty. 
	\label{eq:l_bound}
	\end{equation}
\end{itemize}
It is worth noting that the quantity $R$ will be used only in the convergence analysis, but not in the algorithm itself. As it was pointed in \cite{yurtsever2015universal}, the first assumption is reasonable. Namely, if the set $Q$ is bounded, then $\nabla \vp(\lambda)$ is bounded and H\"older condition holds with $\nu = 0$. If $f(x)$ is uniformly convex, i.e., for all $x,y \in Q$, $\la \nabla f(x) - \nabla f(y) \ra \geqslant \mu \|x-y\|^{\rho}$, for some $\mu > 0$, $\rho \geqslant 2$, then $\nabla \vp(\lambda)$ satisfies H\"older condition with $\nu = \frac{1}{\rho-1}$, $M_{\nu} = \left(\frac{\|\bm{A}\|_{E \to H}^2}{\mu}\right)^{\frac{1}{\rho-1}}$. Here the norm of an operator $\bm{A}:E_1 \to E_2$ is defined as follows
$$
\|\bm{A}\|_{E_1 \to E_2} = \max_{x \in E_1,u \in E_2^*} \{\la u, \bm{A} x \ra : \|x\|_{E_1} = 1, \|u\|_{E_2,*} = 1 \}.
$$

We choose Euclidean proximal setup in the dual space, which means that we introduce Euclidean norm $\|\cdot \|_2$ in the space of vectors $\lambda$ and choose the prox-function $d(\lambda) = \frac12\|\lambda\|_2^2$. Then, we have for the Bregman distance $V[\zeta](\lambda) = \frac12\|\lambda-\zeta\|_2^2$.
Our primal-dual algorithm for Problem $(P_1)$ is listed below as Algorithm \ref{Alg:PDULSGD}. 

\begin{algorithm}[h!]
\caption{PDUGDsDR}
\label{Alg:PDULSGD}
{\small
\begin{algorithmic}[1]
   \REQUIRE starting point $\lambda_0 = 0$, accuracy $\tilde{\eps}_f,\tilde{\eps}_{eq} > 0$.
   \STATE Set $k=0$, $A_0=\alpha_0=0$, $\eta_0=\zeta_0=\lambda_0=0$.
   \REPEAT
		\STATE $\beta_k = \argmin_{\beta \in \left[0, 1 \right]} \vp\left(\zeta^k + \beta (\eta^k - \zeta^k)\right)$; $\lambda^k = \zeta^k + \beta_k (\eta^k - \zeta^k)$ 
		\STATE $h_{k+1} = \argmin_{h \geqslant 0} \vp\left(\lambda^k - h\nabla \vp(\lambda^k)\right)$; $\eta^{k+1} = \lambda^k- h_{k+1}\nabla \vp(\lambda^k)$ // Choose $\nabla \vp(\lambda^k)$ : $\langle \nabla \vp(\lambda^k), \zeta^k - \lambda^k \rangle \geqslant 0$
		\STATE Choose $a_{k+1}$ from $\vp(\eta^{k+1}) = \vp(\lambda^k) - \frac{a_{k+1}^2}{2A_{k+1}} \|\nabla \vp(\lambda^k) \|_2^2 + \frac{\varepsilon a_{k+1}}{2A_{k+1}}$ // $A_{k+1} = A_{k} + a_{k+1}$
		\STATE $\zeta^{k+1} = \zeta^k - a_{k+1}\nabla \vp(\lambda^k)$
		\STATE Set
				\begin{equation}
					\hat{x}^{k+1} = \frac{1}{A_{k+1}}\sum_{i=0}^{k} a_{i+1} x(\lambda^i) = \frac{a_{k+1}x(\lambda^{k})+A_k\hat{x}^{k}}{A_{k+1}}.
				\notag
				\end{equation}
			\STATE Set $k=k+1$.
  \UNTIL{$|f(\hat{x}^{k+1})+\vp(\eta^{k+1})| \leqslant \tilde{\eps}_f$, $\|\bm{A}\hat{x}^{k+1}-b\|_{2} \leqslant \tilde{\eps}_{eq}$.}
	\ENSURE The points $\hat{x}^{k+1}$, $\eta^{k+1}$.	
\end{algorithmic}
}
\end{algorithm}

\begin{Th}
\label{Th:PD_rate}
Let the objective $\vp$ in the problem $(P_2)$ have H\"older-continuous subgradient and the solution of this problem be bounded, i.e. $\|\lambda^*\|_2 \leqslant R$. Then, for the sequence $\hat{x}^{k+1},\eta^{k+1}$, $k\geqslant 0$, generated by Algorithm \ref{Alg:PDULSGD}, 
\begin{align}
&\|\bm{A} \hat{x}^k - b \|_2 \leqslant \frac{2R}{A_k}+ \frac{\eps}{2R}, \quad |\vp(\eta^k) + f(\hat{x}^k)| \leqslant \frac{2R^2}{A_k}+ \frac{\eps}{2},
\label{eq:untileq}
\end{align}
where $A_k \geqslant \left[\frac{1+\nu}{1-\nu}\right]^\frac{1-\nu}{1+\nu}\frac{k^\frac{1+3\nu}{1+\nu}\eps^\frac{1-\nu}{1+\nu}}{2^\frac{1+3\nu}{1+\nu}M_\nu^\frac{2}{1+\nu}}$.
\label{ThLS:PDASTMConv}
\end{Th}
Let us make a remark on complexity. As it can be seen from Theorem \ref{Th:PD_rate}, whenever $A_k \geqslant 2R^2/\eps$, the error in the objective value and equality constraints is smaller than $\eps$. At the same time, using the lower bound for $A_k$, we obtain that the number of iterations to achieve this accuracy is $O\left(\left(\frac{M_\nu^\frac{2}{1+\nu}R^2}{\eps^\frac{2}{1+\nu}} \right)^{\frac{1+\nu}{1+3\nu}}\right)$. Since the algorithm does not use the value of $\nu$, we can take infimum in $\nu \in [0,1]$ of this complexity. This means that the method is uniformly optimal for the class of problems with H\"older-continuous gradient.

\section{Conclusion}

This thesis is based on published papers \cite{dvurechensky2016stochastic,gasnikov2016stochasticInter,bogolubsky2016learning,dvurechensky2020accelerated,dvurechensky2015primal-dual,chernov2016fast,dvurechensky2018computational,dvurechensky2018decentralize,guminov2019accelerated,nesterov2020primal-dual}. 

In papers \cite{dvurechensky2016stochastic,gasnikov2016stochasticInter,bogolubsky2016learning,dvurechensky2020accelerated} we developed optimization methods with (stochastic) inexact first-order oracle, inexact zero-order oracle, inexact directional derivative oracle. We also considered a particular application to learning a parametric model for web-page ranking.

Papers \cite{dvurechensky2015primal-dual,chernov2016fast,dvurechensky2018computational,dvurechensky2018decentralize,guminov2019accelerated,nesterov2020primal-dual} devoted to primal dual methods for convex problems with linear constraints. In particular, we consider infinite-dimensional problems and propose dimension-independent convergence rates for this problem. We also consider (stochastic) convex problems with linear constraints and propose accelerated gradient methods with optimal convergence rates. We apply these methods for approximating optimal transport distance and barycenters. 

Let us list the main results that are obtained in this thesis and submitted for defense.
\begin{enumerate}
    \item Stochastic intermediate gradient method for convex problems with stochastic inexact oracle.
    \item Gradient method with inexact oracle for deterministic non-convex optimization and gradient-free method with inexact oracle for deterministic convex optimization.
    \item A concept of inexact oracle for the methods which use directional derivatives, accelerated and non-accelerated inexact directional derivative method for strongly convex smooth stochastic optimization.
    \item Primal-dual methods for solving infinite-dimensional games in convex-concave and strongly convex-concave setting.
    \item Non-adaptive and adaptive accelerated primal-dual gradient method for strongly convex minimization problems with linear equality and inequality constraints.
    \item New complexity estimates for the optimal transport problem.
    \item Stochastic primal-dual accelerated gradient method for problems with linear constraints and its application to the problem of approximation of Wasserstein barycenter.
    \item A universal primal-dual accelerated gradient method with line-search.
\end{enumerate}

\section*{Acknowledgements}
The dissertation was supported by Russian Science Foundation (project 18-71-10108), by RFBR project number 18-31-20005 mol-a-ved, by the Ministry of Science and Higher Education of the Russian Federation (Goszadaniye) No.075-00337-20-03, project No. 0714-2020-0005, and by RFBR project number 18-29-03071-mk.

\bibliographystyle{plainnat}
\bibliography{PD_references,my_papers,main}

\iftoggle{thesis}{

\newpage

\section{Appendix}

\appendix

\section{Статья ``Fast gradient descent for convex mi\-ni\-mi\-zation problems with an oracle producing a $(\delta, L)$-model of function at the requested point''}

\hspace{1cm}

\textbf{Авторы:} Gasnikov A., Tyurin A.

\textbf{Аннотация:} A new concept of $(\delta, L)$--model of a function that is a generalization of the Devolder–Glineur–Nesterov $(\delta, L)$--oracle is proposed. Within this concept, the gradient descent and fast gradient descent methods are constructed and it is shown that constructs of many known methods (composite methods, level methods, conditional gradient and proxi\-mal meth\-ods) are particular cases of the methods proposed in this paper.

\includepdf[pages=-]{1_2016_Dvurechensky-Gasnikov_Stochastic_Intermediate_Gradient_JOTA.pdf}

\section{Статья ``Gradient methods for problems with in\-exact model of the objective''}

\hspace{1cm}

\textbf{Авторы:} Stonyakin F., Dvinskikh D., Dvurechensky P., Kroshnin A., Kuznetsova O., Agafonov A., Gasnikov A., Tyurin A., Uribe С., Pasechnyuk D., Artamonov S.

\textbf{Аннотация:} We consider optimization methods for convex minimization problems under inexact information on the objective function. We introduce inexact model of the objective, which as a particular cases includes inexact oracle [16] and relative smoothness condition [36]. We analyze gradient meth\-od which uses this inexact model and obtain convergence rates for convex and strongly convex problems. To show potential applications of our general framework we consider three particular problems. The first one is clustering by electorial model introduced in [41]. The second one is approximating optimal transport distance, for which we propose a Proximal Sinkhorn algo\-rithm. The third one is devoted to approximating optimal transport barycen\-ter and we propose a Proximal Iterative Bregman Projections algorithm. We also illustrate the practical performance of our algorithms by numerical experiments.

\includepdf[pages=-]{2_2016_Gasnikov-Dvurechensky_StochasticIntermediateGradient_Doklady.pdf}

\section{Статья ``Primal-dual fast gradient method with a model''}

\hspace{1cm}

\textbf{Авторы:} Tyurin A.

\textbf{Аннотация:} In this work we consider a possibility to use the conception of $(\delta, L)$--model of a function for optimization tasks, whereby solving a primal problem there is a necessity to recover a solution of a dual problem. The conception of $(\delta, L)$--model is based on the conception of $(\delta, L)$--oracle which was proposed by Devolder–Glineur–Nesterov, herewith the authors proposed approximate a function with an upper bound using a convex quadratic function with some additive noise $\delta$. They managed to get convex quadratic upper bounds with noise even for nonsmooth functions. The conception of $(\delta, L)$--model continues this idea by using instead of a convex quadratic function a more complex convex function in an upper bound. Possibility to recover the solution of a dual problem gives great benefits in different problems, for instance, in some cases, it is faster to find a solution in a primal problem than in a dual problem. Note that primal-dual methods are well studied, but usually each class of optimization problems has its own primal-dual method. Our goal is to develop a method which can find solutions in different classes of optimization problems. This is realized through the use of the conception of $(\delta, L)$--model and adaptive structure of our methods. Thereby, we developed primal-dual adaptive gradient method and fast gradi\-ent method with $(\delta, L)$--model and proved convergence rates of the methods, moreover, for some classes of optimization problems the rates are optimal. The main idea is the following: we find a dual solution to an approximation of a primal problem using the conception of $(\delta, L)$--model. It is much easier to find a solution to an approximated problem, however, we have to do it in each step of our method, thereby the principle of “divide and conquer” is realized.

\includepdf[pages=-]{3_2016_Bogolubski-learning-supervised-pagerank-with-gradient-based-and-gradient-free-optimization-methods_supplementary.pdf}

\section{Статья ``Accelerated and nonaccelerated stocha\-stic gradient descent with model conception''}

\hspace{1cm}

\textbf{Авторы:} Dvinskikh D., Tyurin A., Gasnikov A., Omelchenko S.

\textbf{Аннотация:} In this paper, we describe a new way to get convergence rates for optimal methods in smooth (strongly) convex optimization tasks. Our approach is based on results for tasks where gradients have nonrandom small noises. Unlike previous results, we obtain convergence rates with model conception.

\textbf{Примечание:} статья принята в печать в журнал Mathematical Notes V. 108, № 4. Ниже представлен препринт arXiv: 2001.03443.

\includepdf[pages=-]{4_2020_Dvurechensky_et_al_Accelerated Directional Derivetive.pdf}

\section{Статья ``Heuristic adaptive fast gradient method in stochastic opti\-mization tasks''}

\hspace{1cm}

\textbf{Авторы:} Ogaltsov A., Tyurin A.

\textbf{Аннотация:} A heuristic adaptive fast stochastic gradient descent method is proposed. It is shown that this algorithm has a higher convergence rate in practical problems than currently popular optimization methods. Furthermore, a justification of this method is given, and difficulties that prevent obtaining optimal estimates for the proposed algorithm are described.

\textbf{Примечание:} ниже представлена русская версия статьи в Журнал вычислительной математики и математической физики, 2020, том 60, № 7, с. 1143–1150.

\includepdf[pages=-]{5_2015_Dvurechensky_et al_Primal-DualMethodsForSolvingInfiniteDimensionalGames.pdf}

\section{Статья ``A stable alternative to Sinkhorn’s algo\-rithm for regularized optimal transport''}

\hspace{1cm}

\textbf{Авторы:} Dvurechensky P., Gasnikov A., Omelchenko A., Tyurin A.

\textbf{Аннотация:} In this paper, we are motivated by two important applica\-tions: entropy-regularized optimal transport problem and road or IP traffic demand matrix estimation by entropy model. Both of them include solving a special type of optimization problem with linear equality constraints and objective given as a sum of an entropy regularizer and a linear function. It is known that the state-of-the-art solvers for this problem, which are based on Sinkhorn’s method (also known as RSA or balancing method), can fail to work, when the entropy-regularization parameter is small. We consider the above optimization problem as a particular instance of a general strongly convex optimization problem with linear constraints. We propose a new algorithm to solve this general class of problems. Our approach is based on the transition to the dual problem. First, we introduce a new accelerated gradient method with adaptive choice of gradient’s Lipschitz constant. Then, we apply this method to the dual problem and show, how to reconstruct an approximate solution to the primal problem with provable convergence rate. We prove the rate $O(1/k^2)$, k being the iteration counter, both for the absolute value of the primal objective residual and constraints infeasibility. Our method has similar to Sinkhorn’s method complexity of each iteration, but is faster and more stable numerically, when the regularization parameter is small. We illustrate the advantage of our method by numerical experiments for the two mentioned applications. We show that there exists a threshold, such that, when the regularization parameter is smaller than this threshold, our method outper\-forms the Sinkhorn’s method in terms of computation time.

 \includepdf[pages=-]{6_2016_Chernov_Fast primal-dual gradient method for strongly convex minimization problems with linear constraints_chapter.pdf}

\section{Статья ``Dual approaches to the minimization of strongly convex functionals with a simple struc\-ture under affine con\-stra\-ints''}

\hspace{1cm}

\textbf{Авторы:} Anikin A.,  Gasnikov A., Dvurechensky P., Tyurin A., Chernov A.

\textbf{Аннотация:} A strongly convex function of simple structure (for example, separable) is minimized under affine constraints. A dual problem is construct\-ed and solved by applying a fast gradient method. The necessary properties of this method are established relying on which, under rather general conditions, the solution of the primal problem can be recovered with the same accuracy as the dual solution from the sequence generated by this method in the dual space of the problem. Although this approach seems natural, some previously unpublished rather subtle results necessary for its rigorous and complete theoretical substantiation in the required generality are presented.

\includepdf[pages=-]{7_2018_Dvurechensky_et al_Computational optimal transport Complexity by accelerated gradient descent.pdf}
\includepdf[pages=-]{7_2018_Dvurechensky_et al_Computational optimal transport Complexity by accelerated gradient descent_supp.pdf}

\section{Статья ``Accelerated gradient sliding for minimizing the sum of functions''}

\hspace{1cm}

\textbf{Авторы:} Dvinskikh D., Omelchenko A., Gasnikov A., Tyurin A.

\textbf{Аннотация:} In this article, we propose a new way to justify the accelerated gradient sliding of G. Lan, which allows one to extend the sliding technique to a combination of an accelerated gradient method with an accelerated variance reduced method. We obtain new optimal estimates for solving the problem of minimizing a sum of smooth strongly convex functions with a smooth regularizer. 

 \includepdf[pages=-]{8_2018_Dvurechensky_et al_Decentralize-and-randomize-faster-algorithm-for-wasserstein-barycenters_full.pdf}

 \section{Статья ``Accelerated gradient sliding for minimizing the sum of functions''}

\hspace{1cm}

\textbf{Авторы:} Dvinskikh D., Omelchenko A., Gasnikov A., Tyurin A.

\textbf{Аннотация:} In this article, we propose a new way to justify the accelerated gradient sliding of G. Lan, which allows one to extend the sliding technique to a combination of an accelerated gradient method with an accelerated variance reduced method. We obtain new optimal estimates for solving the problem of minimizing a sum of smooth strongly convex functions with a smooth regularizer. 

 \includepdf[pages=-]{9_2019_Guminov_Accelerated.pdf}

\section{Статья ``Accelerated gradient sliding for minimizing the sum of functions''}

\hspace{1cm}

\textbf{Авторы:} Dvinskikh D., Omelchenko A., Gasnikov A., Tyurin A.

\textbf{Аннотация:} In this article, we propose a new way to justify the accelerated gradient sliding of G. Lan, which allows one to extend the sliding technique to a combination of an accelerated gradient method with an accelerated variance reduced method. We obtain new optimal estimates for solving the problem of minimizing a sum of smooth strongly convex functions with a smooth regularizer. 

 \includepdf[pages=-]{10_2020_Nesterov_et_al_Primal_dual.pdf}

}  


\end{document}